\documentclass[12pt,leqno,color,amssymb]{amsart}
\usepackage{amsthm,amsmath,amssymb,latexsym,graphics}
\usepackage{graphicx}
\usepackage{enumerate}
\usepackage{psfrag}
\newcommand{\field}[1]{\mathbb{#1}}\newcommand{\R}{\field{R}}

\newcommand{\Hp}{\field{H}}
\newcommand{\N}{\field{N}}
\newcommand{\C}{\field{C}}
\newcommand{\Q}{\field{Q}}
\newcommand{\Z}{\field{Z}}

\newcommand{\ML}{\mathcal{ML}}

\newcommand{\eps}{\varepsilon}
\newcommand{\CB}{\mathcal{B}}
\newcommand{\CN}{\mathcal{N}}

\newcommand{\GF}{\mathcal{GF}}
\newcommand{\CR}{\mathcal{R}}
\newcommand{\IM}{\mathring{M}}
\theoremstyle{plain}
\newtheorem{theorem}{Theorem}[section]

\newtheorem{prop}[theorem]{Proposition}
\newtheorem{lemma}[theorem]{Lemma}

\newtheorem{claim}[theorem]{Claim}
\newenvironment{rem}{\bf Remark.\rm}{\hfill $\diamond$\\}

\theoremstyle{definition}
\newtheorem{note}[theorem]{Note}
\newtheorem{definition}[theorem]{Definition}
\theoremstyle{plain}
\begin{document}
\title{Properness of the bending map}
\author{Cyril Lecuire}
\address{\hskip-\parindent
	UMPA UMR 5669 CNRS\\  	  	 
	ENS de Lyon Site Monod\\ 	  	  	 
	46 Allée d'Italie\\  	 
	69364 Lyon Cedex 07\\
	France}
\email{cyril.lecuire@ens-lyon.fr}

\maketitle


\begin{abstract}
The bending map of a hyperbolic $3$-manifold with boundary  maps a geometrically hyperbolic metric to its bending measured geodesic lamination. We show that the bending map is proper. As a byproduct of the proof we show that the group of isotopy classes of homeomorphisms of $M$ acts properly discontinuously on the set of doubly incompressible measured geodesic laminations.
\end{abstract}

\section{Introduction}
Consider a compact, orientable $3$-manifold $M$ that is hyperbolic, namely the interior $\mathring{M}$ of $M$ is endowed with a complete metric $\sigma$ of constant sectional curvature $-1$. A fundamental subset of $(M,\sigma)$ is its convex core $N(\sigma)$, defined as the smallest non-empty closed subset of the interior of $M$ that is locally convex (with respect to $\sigma$) and homotopy equivalent to $M$. Its boundary $\partial N(\sigma)$, when not empty, is almost everywhere totally geodesic and is bent along a geodesic lamination. The amount of bending is described by a measured geodesic lamination called the {\underline\it bending measured geodesic lamination} of $\sigma$ (cf. \cite{thurston:notes} or \cite{dbaepstein:marden}). This yields a bending map $b:\mathcal{GF}(M)\rightarrow \mathcal{ML}(\partial M)$ which to a hyperbolic metric associates its bending measured geodesic lamination. In the present paper, we investigate the properness properties of this bending map.

A complete hyperbolic metric $\sigma$ on $\IM$ is {\em geometrically finite} if $N(\sigma)$ has finite volume. We denote by $\GF(M)$ the set of isotopy classes of geometrically finite metrics on $\IM$, two metrics $\sigma_1,\sigma_2$ on $\IM$ are isotopic if there is a diffeomorphism  $\varphi:\IM\to\IM$ isotopic to the identity such that $\sigma_2=\varphi_*\sigma_1$. We exclude from $\GF(M)$ the Fuchsian metrics, i.e. the metrics with two-dimensional convex cores. We equip $\mathcal{GF}(M)$ with the marked pointed Hausdorff-Gromov topology. Let us choose a point $x$ in $int(M)$. A metric $\sigma_2$ lies in a $(k,r)$-neighborhood of $\sigma_1$ if there exists a diffeomorphism $\phi:M\rightarrow M$ isotopic to the identity such that the restriction of $\phi$ to the ball $B(x,r)\subset (M,\sigma_1)$ is a $(k,\eps)$-quasi-isometry into its image in $(M,\sigma_2)$. We obtain a basis of neighborhoods of $\sigma_1$ by letting $k$ tend to $1$, $\eps$ tend to $0$ and $r$ tend to $+\infty$. The topology defined in this way does not depend on the choice of the point $x$ and if we replace $(k,\eps)$-quasi-isometry by $k$-biLipschitz map, we obtain the same topology (see \cite[\textsection E.1]{benedetti:petronio}).

We have already introduced the bending map $b:\mathcal{GF}(M)\rightarrow \mathcal{ML}(\partial M)$. Its image, $b(\mathcal{GF}(M))$ has been described in \cite{bonahon:otal} and \cite{lecuire:plissage}, it is the set $\mathcal{P}(M)$ of measured geodesic laminations satisfying the following conditions:
\begin{enumerate}[a)]
\item the weigth of any closed leaf of $\lambda$ is at most $\pi$;
\item $\exists\eta>0$ such that, for any essential annulus $E$, $i(\partial E,\lambda)\geq\eta$;
\item $i(\lambda,\partial D)>2\pi$ for any essential disc $D$.
\end{enumerate}

The space of measured geodesic lamination is usually equipped with the weak$^*$ topology (see \cite[\textsection A.3]{otal:fibre}) and this provides a natural topology on $\mathcal{P}(M)$. The continuity and the differentiability of $b_\mathcal{CC}$ have been studied in \cite{keen:series:continuity}, \cite{lecuire:continuity} and \cite{bonahon:variations}. Its injectivity has been established in \cite{bonahon:otal} in the special case of weighted multi-curves, in \cite{bonahon:almost:fuchsian} near the Fuchsian locus and in \cite{series:bending:conjecture} for punctured torus groups and very recently in the case of convex cocompact groups in \cite{dular:schlenker:preprint}. The bending map naturally extends to a map mixing bending laminations and ending laminations and the properness and the image of this extension have been studied in \cite{baba:ohshika} for Kleinian surface groups. In the present paper we will study the properness of $b_\mathcal{GF}$.

A complete hyperbolic metric $\sigma$ on $\IM$ is {\em convex cocompact} if $N(\sigma)$ is compact. Let $\mathcal{CC}(M)\subset\GF(M)$ denote the set of isotopy classes of convex cocompact metrics. It follows from the definition of the bending measure that a lamination $\lambda\subset \mathcal{P}(M)=b(\mathcal{GF}(M))$ lies in $b(\mathcal{CC}(M))$ if and only it does not have any leaf with a weight equal to $\pi$. Let us state a simpler version of our main theorem:
\begin{theorem}
	The map $b_\mathcal{CC}$ from $\mathcal{CC}(M)$ to $b(\mathcal{CC}(M))$ is proper.
\end{theorem}

As observed in \cite{lecuire:continuity}, when dealing with geometrically finite metrics, we need to consider limiting lamination with weights greater than $\pi$. For this purpose, we define the clipping map $c:\ML(\partial M)\to\ML(\partial M)$ by capping the weights at $\pi$. Namely, the measured geodesic lamination $c(\lambda)$ is obtained by replacing  by $\pi$ the weights of the leaves of $\lambda$ which have a weight greater than $\pi$. We define a relation $\mathcal{R}$ on $\ML(\partial M)$ by  $\lambda \mathcal{R} \mu$ if and only if $c(\lambda)=c(\mu)$. We denote by $\dot\lambda$ the class of $\lambda$ modulo $\CR$. To define the topology of $\ML(\partial M)/\mathcal R$, let us denote by $\ML^+(\partial M)$ the space of extended measured geodesic laminations in which we allow closed leaves to have infinite weights and let us equip $\ML^+(\partial M)$ with a topology that extends the weak$^*$ topology on $\ML(S)$ (see details in \textsection \ref{topologies}). There is a natural bijection between the quotients $\ML(\partial M)/\mathcal R$ and $\ML^+(\partial M)/\mathcal R$ but they inherit different topologies from the quotient maps (see \textsection \ref{properties of MLR}). We call the topology inherited from $\ML^+(\partial M)$ tubular topology. Unless stated otherwise $\ML(\partial M)/\mathcal R$ will be equipped with the tubular topology. 
We define a map $b_\mathcal{R}: \mathcal{GF}(M)\to \mathcal{ML}(\partial M)/\mathcal{R}$ by taking the projection of the image of the bending map: if $b(\sigma)=\lambda$, then $b_{\mathcal R}(\sigma)=\dot\lambda$. Since, as mentioned above, $b(\mathcal{GF}(M))=\mathcal{P}(M)$, the image of $b_{\CR}$ is the projection of $\mathcal P (M)$ to $\mathcal{ML}(\partial M)/\mathcal{R }$ that we will denote by $\mathcal{P}(M)/\mathcal{R}$. Notice that although the restriction of the projection to $\mathcal{P}(M)$ is a bijection, the topology of $\mathcal{P}(M)/\mathcal{R}$ as a subset of $\mathcal{ML}(\partial M)/\mathcal{R }$ equipped with the tubular topology is different from the topology of $\mathcal{P}(M)$. This is illustrated by the fact that the map $b_{\mathcal R}$ is continuous while $b$ is not (see \cite{lecuire:continuity}). Notice that in \cite{lecuire:continuity}, $\mathcal{ML}(\partial M)/\mathcal{R }$ is equipped with the quotient topology. The main Theorem of \cite{lecuire:continuity} is false with this assumption but when we replace the quotient topology with the tubular topology, then the proofs of the main statement and all intermediate statements are correct. In the present article, we complete the description of the behavior of the bending map when one allows new parabolics to appear by proving:

\begin{theorem}		\label{two}
The map $b_\mathcal{R}$ from $\mathcal{GF}(M)$ to $\mathcal{P}(M)/\mathcal{R}$ equipped with the tubular topology is proper.
\end{theorem}


The set $\mathcal{P}(M)$ is closely related to the set of doubly incompressible laminations $\mathcal{D}(M)$ introduced in \cite{lecuire:masur}. A measured geodesic lamination $\lambda\in\mathcal{ML}(\partial M)$ is {\em doubly incompressible} if  $\exists\eta>0$ such that $i(\partial E,\lambda)>\eta$ for any essential annulus or disc  $E$. Some of the arguments used to prove Theorem \ref{two} may also be used to study the action of the modular group $Mod (M)$, the group of isotopy classes of diffeomorphisms of $M$, leading to:

\begin{theorem}		\label{tree}
When $M$ is not a genus two handlebody, the action of $Mod (M)$ on $\mathcal{D}(M)$ is properly discontinuous.
\end{theorem}

This result is also essentially true for a genus two handlebody, it simply needs a little tweaking of the definitions (see condition (-) in \cite{lecuire:masur}).

\indent The paper is organized as follows. In section 2, we state some definitions and facts about geodesic laminations, $\ML(S)/\mathcal R$, convergence of representations, $\mathcal{P}(M)$ and $\mathcal{D}(M)$. In section 3, we explain how to use a bound on the lengths of the bending laminations to get a converging subsequence of representations associated to sequence of metrics. In section 4, we show how to compare the induced metric on the boundary with the metric outside, and to upgrade the convergence of associated representations to something closer to the convergence of metrics. In section 5, we conclude the proof of Theorem \ref{two} and prove Theorem \ref{tree}.

\indent
I would like to thank Francis Bonahon for discussions from which this paper originates, Jean-Pierre Otal for his advices and Young Eun Choi for useful comments. I also thank the anonymous referee whose comments lead to finding an error in \cite{lecuire:continuity} and defining the tubular topology on $\ML(S)/\mathcal R$.
\section{Definitions}

\subsection{Geodesic laminations}
A {\em geodesic lamination} on a finite type hyperbolic surface $S$ is a compact subset which is a disjoint union of complete simple geodesics embedded in $S$. The space of geodesic laminations on $S$ endowed with the Hausdorff topology is denoted by $\mathcal{L}(S)$. A geodesic lamination whose leaves are all closed is called a {\em multi-curve}. If each half-leaf of a geodesic lamination $L$ is dense in $L$, then $L$ is {\em minimal}. Such a minimal geodesic lamination is either a simple closed curve or an {\em irrational lamination}. A leaf $l$ of a geodesic lamination $L$ is {\em recurrent} if it lies in a minimal geodesic sublamination. Any geodesic lamination is the disjoint union of  finitely many minimal laminations and non-recurrent leaves.

It is a classical fact that the definition of $\mathcal{L}(S)$ can be made independent of the choice of the hyperbolic metric on $S$. This is explained in the following note (compare with \cite[\textsection I.4.1.4]{canary:epstein:green}).

\begin{note}
Consider two complete hyperbolic metrics $s_1$ and $s_2$ on a finite type surfacet $S$ and let $q:\Hp^2\rightarrow S $ be a covering map. Let $l\subset S $ be a $s_1$-geodesic and let $\hat l\subset q^{-1}(l)$ be a lift of $l$. There is  a unique $s_2$-geodesic $\Theta(l)$ such that there is a bounded homotopy between a lift $\hat\Theta(\hat l)\subset q^{-1}(\Theta(l))$ of $\Theta(l)$ and $\hat l$. Here a bounded homotopy is a map $F:\R\times[0,1]\rightarrow S$ such that the lengths of the arcs $F(\{x\}\times[0,1])$ are uniformly bounded. So we get a homeomorphism $\Theta :\{\mbox{geodesics of }(S ,s_1)\}\rightarrow\{\mbox{geodesics of }(S,s_2)\}$ and the image $\Theta(L )$ of a $s_1$-geodesic lamination $L $ is a $s_2$-geodesic lamination. We will say that $\Theta(L )$ is the geodesic lamination $L $ for the metric $s_2$. Thus the geodesic lamination $L $ is well defined for any hyperbolic metric on $S$.\\
\indent
Let $l_1 $ and $l_2 $ be two $s_1$-geodesic laminations and let $x\subset l_1 \cap l_2 $ be a transverse intersection. Let $\hat x\subset\Hp^2$ be a lift of $x$. This point $\hat x$ is the intersection of a leaf $\hat l_1$ of $q^{-1}(l_1 )$ and of a leaf $l_2$ of $q^{-1}(l_2 )$. Since $\hat\Theta(\hat l_1)$ intersects  $\hat\Theta(\hat l_2)$ transversely, there is a map $l_1 \cap l_2\rightarrow \Theta(l_1) \cap \Theta(l_2)$ which to $x$ associates $q(\hat\Theta(\hat l_1)\cap\hat\Theta(\hat l_2))$. This way, each point lying in the transverse intersection of two geodesic laminations is defined independently of the choice of the hyperbolic metric on $S$.
\end{note}

Given a connected geodesic lamination $L$ on a surface $S$ of finite type, the {\em minimal supporting surface of $L$}, $S(L)$ is the smallest essential compact subsurface containing $L$. By {\em essential}, we mean that no component of $\partial S(L)$ bounds a disc in $S$. This surface $S(L)$ is well defined up to isotopy. Notice that two different component of $\partial S(L)$ may bound an annulus. In particular, when $S$ is endowed with a complete hyperbolic metric, we cannot expect $\partial S(L)$ to be geodesic. On the other hand each component of $\partial S(L)$ is isotopic to a simple closed geodesic and we will also denote by $\partial S(L)$ the collection of those geodesics. When $L$ is not connected, we take $S(L)$ to be the disjoint union of the minimal supporting surfaces of the connected components of $L$.

\subsubsection{Measured laminations}

A  {\em measured geodesic lamination} $\lambda$ is a transverse measure for some geodesic lamination $|\lambda|$. Any arc $k\approx [0,1]$ embedded in $S$ transversely to $|\lambda|$, such that $\partial k\subset S\setminus\lambda$, is endowed with an  additive measure $d\lambda$ such that:\\
\indent - the support of $d\lambda_{|k}$ is $|\lambda|\cap k$;\\
\indent - if an arc $k$ can be homotoped into $k'$ by a homotopy respecting $|\lambda|$ then $\int_k\! d\lambda=\int_{k'} d\lambda$.

An arc $k\subset S$ is generic if any intersection between $k$ and a simple geodesic is transverse. It follows from \cite{birman:series} that any arc can be approximated by generic arcs. We denote by $\mathcal{ML}(S)$ the space of measured geodesic lamination equipped with the weak$^*$ topology that can be defined by the family of semi-norms $\lambda\to|\int_k d\lambda|$ where $k$ ranges over all generic arcs. A sequence $(\lambda_n)\in\ML(S)$ converges to $\lambda\in\ML(S)$ if and only if $\int_k\lambda_n\longrightarrow\int_k\lambda$ for every generic arc $k\subset S$. It has been observed by Thurston that the weak$^*$ topology can be defined using only finitely many arcs, \cite[Theorem 16]{bonahon:laminations}, and that $\ML(S)$ is a finite dimensional piecewise linear manifold, \cite[Proposition 9.5.8]{thurston:notes}.\\

Let $\gamma$ be a weighted simple closed geodesic with support $|\gamma|$ and weight $w$ and let $\lambda$ be a measured geodesic lamination. The intersection number of $\gamma$ and $\lambda$ is defined by $i(\gamma,\lambda)=w \int_{|\gamma|} d\lambda$. If $|\gamma|$ is a leaf of $\lambda$, then we define the intersection number by $i(\gamma,\lambda)=0$. The weighted simple closed curves are dense in $\mathcal{ML}(S)$ (\cite[proposition 8.10.7]{thurston:notes}) and this intersection number extends continuously to a function $i:\mathcal{ML}(S)\times\mathcal{ML}(S)\rightarrow\R$ (cf. \cite{bonahon:bouts}).

\subsubsection{The space $\ML(S)/\mathcal R$}	\label{topologies}

In the introduction, we have defined the
relation $\mathcal{R}$ on $\ML(S)$: $\lambda \mathcal{R} \mu$ if and only if $c(\lambda)=c(\mu)$ where the measured geodesic lamination $c(\lambda)$ (resp. $c(\mu)$) is obtained by replacing  by $\pi$ the weights of the closed leaves of $\lambda$ (resp. $\mu$) which have a weight greater than $\pi$. Given $\lambda\in\ML(S)$, we denote by $\dot\lambda$ its projection in $\ML(S)/\mathcal{R}$. Conversely, we define a map $\breve c:\ML(S)/\mathcal R\to\ML(S)$ as follows: given $\lambda\in\ML(S)$ that projects to $\dot\lambda\in \ML(S)/\mathcal R\to\ML(S)$, $\breve c(\dot\lambda)=c(\lambda)$. In other words $\breve c(\dot\lambda)$ is the representative of the equivalence class $\dot\lambda$ that has no leaf with a weight greater than $\pi$.\\

As explained in the introduction, rather than the quotient topology, we will use what we call the {\em tubular topology} on $\ML(S)/\mathcal R$. To prepare its definition, we are going to build a specific basis of neighborhoods for the weak$^*$ topology on $\ML(S)$ based on a construction due to Otal,\cite[\textsection 3.2]{otal:fibre}. 

Let $\lambda$ be a measured geodesic lamination on a hyperbolic surface $S$ of finite area. We are going to construct a finite family of generic arcs that will define a local basis at $\lambda$. Let $L$ be a connected component of its support $|\lambda|$. If $L$ is a simple closed curve with a weight equal to or larger than $\frac{\pi}{2}$, we choose a generic arc that intersects $L$ once and is disjoint from the other components of $L$. Otherwise, we start with a generic arc $k$ that intersect $L$ such that $\int_{k} d\lambda<\frac{\pi}{2}$. By \cite[Proposition A.3.4]{otal:fibre}, every leaf of $L$ is dense in $L$ and by \cite[\textsection 3.2]{otal:fibre} the closures of the components of $L\setminus k$ form finitely many families such that any two arcs in a family are isotopic relative to $k$. As explained in \cite[\textsection 3.2]{otal:fibre}, this allows us to build finitely many rectangles $r_i:[0,1]\times[0,1]\to S$ such that $r_i(\{\eta\}\times[0,1])\subset k$ for $\eta\in\{0,1\}$, $r_i([0,1]\times\{\eta\})$ is a generic arc for $\eta\in\{0,1\}$, the restriction of $r_i$ to $(0,1)\times[0,1]$ is an embedding, $r_i((0,1)\times[0,1])\cap r_j((0,1)\times[0,1])=\emptyset$ for $i\neq j$ and  $L\subset\bigcup_i r_i([0,1]\times(0,1))$. We get a finite family of generic arcs by taking the sides $r_i(\{\eta\}\times[0,1])$ and $r_i([0,1]\times\{\eta\})$, $\eta\in\{0,1\}$, of all the rectangles $r_i$. We do this construction for each component of $|\lambda|$ and obtain a finite family of generic arcs $k_i$ such that $\int_{k_i} d\lambda<\frac{\pi}{2}$ unless $k_i$ intersects a closed leaf of $\lambda$. We add a generic arc in each component of $S\setminus|\lambda|$ to obtain a finite family of generic arcs $k_i, i\leq q$ such that any simple geodesic intersects at least one arc $k_i$.
Set $\CB_\eps(\lambda)=\{\mu\in\ML(S):|\int_{k_i}d\mu-\int_{k_i} d\lambda|<\eps\; \forall i\leq q\}$. It follows from the work of Thurston \cite[\textsection 8.2]{thurston:notes} that $\{\CB_\eps(\lambda)|\eps>0\}$ is a local basis at $\lambda$ for the weak$^*$ topology.

Given $\dot\lambda\in\ML(S)/\CR$, let $F_{\dot\lambda}$ be the subset of $\ML(S)$ made up of all measured geodesic laminations that project to $\dot\lambda$ and let $\lambda^{(p)}$ be the union of the leaves of $\breve c(\dot\lambda)$ with a weight equal to $\pi$. Using the previous construction, we define the set $\breve\CB_\eps(\dot\lambda)=\{\mu\in\ML(S):\int_{k_i}d\mu>\pi-\eps$ if $k_i$ intersects $\lambda^{(p)}$ and $|\int_{k_i}d\mu-\int_{k_i} d\lambda|<\eps$ otherwise$\}$ which is an open set containing $F_{\dot\lambda}$.  Let $m(\dot\lambda)$ be the maximum of the weights of the closed leaves of $\breve c(\dot\lambda)$ with a weight smaller than $\pi$. For $\eps<\min\{\frac{\pi}{2},\pi-m(\lambda)\}$, if $\mu\in \breve\CB_\eps(\dot\lambda)$ then $\int_{k_i}d\mu<\pi$ for any arc $k_i$ that is disjoint from $\lambda^{(p)}$. It follows that any leaf of a measured geodesic lamination $\mu\in \breve\CB_\eps(\dot\lambda)$ with a weight greater than or equal to $\pi$ is a leaf of $\lambda^{(p)}$. In particular $\breve\CB_\eps(\dot\lambda)$ is saturated with respect to $\CR$. Let $\dot\CB_\eps(\dot\lambda)$ be the projection of $\breve\CB_\eps(\dot\lambda)$ to $\ML(S)/\CR$. We define the tubular topology by setting that $\{\dot\CB_\eps(\dot\lambda): 0<\eps<\min\{\frac{\pi}{2},\pi-m(\lambda)\}\}$ is a local basis at $\dot\lambda$. 
From now on $\ML(S)/\CR$ will be equipped with the tubular topology unless stated otherwise.\\

Notice that we have established the following fact that we shall put into a claim for later use:

\begin{claim}	\label{larger than pi}
Consider $(\lambda_n)\subset\ML(S)$ such that $\dot\lambda_n$ converges to $\dot\lambda$. Then, for $n$ large enough, any closed leaf of $\lambda_n$ with a weight equal to or larger than $\pi$ is a closed leaf of $\breve c(\dot\lambda)$ with a weight equal to $\pi$.
\end{claim}

From the definition and the work of Thurston, \cite[\textsection 8 and 9]{thurston:notes} (see also \cite{penner:harer} and \cite{otal:fibre}) we deduce the following characterization of converging sequences in the tubular topology. Consider a sequence $(\dot\lambda_n)\subset \ML(S)/\CR$ and an element $\dot\lambda\in\ML(S)/\CR$ and denote by $\lambda^{(p)}$ the union of the closed leaves of $\breve c(\dot\lambda)$ with weights equal to $\pi$. The sequence $(\dot\lambda_n)$ converges to $\dot\lambda$ in the tubular topology if and only if $(\breve c(\dot\lambda_n))$ converges to $\breve c(\dot\lambda)$ when restricted to $S\setminus \lambda^{(p)}$, i.e. $\int_k d(\breve c(\dot\lambda_n))\longrightarrow\int_k d(\breve c(\dot\lambda))$ for any generic arc $k\subset  S\setminus \lambda^{(p)}$, and $\liminf \int_k d(\breve c(\dot\lambda_n))\geq\pi$ for any generic arc $k$ that intersects $\lambda^{(p)}$.\\

In the introduction we used a shortcut to define this topology by using the space $\ML^+(S)$ of extended measured geodesic laminations where we allow closed leaves to have an infinite weight. Let us explain how this gives an alternate definition of the tubular topology. An element $\lambda$ of $\ML^+(S)$ is a transverse measure with values in $\R^+\cup\{\infty\}$ supported by a geodesic lamination $|\lambda|$ which only gives infinite measures to arcs that intersect some closed leaves (the leaves with infinite weight). We also require $|\lambda|$ to be the support of an element of $\ML(S)$ to avoid having leaves spiraling around a closed leaf with infinite weight. Given $\lambda\in\ML^+(S)$, let us denote by $\lambda^\infty$ the union of its closed leaves with infinite weight and set $\CB_{\eps,K}(\lambda)=\{\mu\in\ML^+(S):\int_{k_i}d\mu\geq K$ if $k_i$ intersects $\lambda^\infty$ and $|\int_{k_i}d\mu-\int_{k_i}d\lambda|<\eps$ otherwise$\}$ where $\{k_i, i\leq q\}$ is the family of arcs constructed above. We extend the weak$^*$ topology to $\ML^+(S)$ by setting that $\{\CB_{\eps,K}(\lambda)|\eps>0,K>0\}$ is a local basis at $\lambda$. It is not difficult to see that the quotient topology on $\ML^+(S)/\CR$ is equivalent to the tubular topology defined above on $\ML(S)/\CR$.

\subsubsection{Properties of $\ML(S)/\mathcal R$}	\label{properties of MLR}

To illustrate the relevance of the tubular topology, let us show that it differs from the quotient topology on $\ML(S)/\CR$. Let $\gamma\in\ML(S)$ be a weighted simple closed curve with weight $\pi$. Let $\{k_i,i\leq q\}$ be the family of arcs constructed in the previous subsection, ordered so that $k_1$ intersects $|\gamma|$. The set $\mathcal V=\{\mu\in \ML(S):\int_{k_2} d\mu<(\int_{k_1} d\mu)^{-1}\}$ is an open subset of $\ML(S)$ that contains $F_{\dot\gamma}$ but does not contain $\breve\CB_\eps(\dot\gamma)$ for any $\eps>0$. For any $\eps<\frac{\pi}{2}$, $\mathcal V\cap \breve\CB_\eps(\dot\gamma)$ is saturated with respect to $\mathcal R$. Hence it projects to $\ML(S)/\CR$ to a subset that is open for the quotient topology but not for the tubular topology. It follows that the quotient topology is finer than the tubular topology.\\

By definition, $\dot{\CB}_{\frac{1}{n}}(\dot\lambda)$ is a countable local basis at $\dot\lambda$. Thus $\ML(S)/\mathcal R$ is first countable. In particular its topology is determined by its converging sequences. With that in mind we are going to establish some separation properties of $\ML(S)/\CR$. 

\begin{claim}
The space $\ML(S)/\CR$ equipped with the tubular topology is regular.
\end{claim}

\begin{proof}
Let us first show by contradiction that it is Hausdorff. Consider $\dot\lambda\neq\dot\mu\in\ML(S)/\CR$, and assume that $\dot{\CB}_{\frac{1}{n}}(\dot\lambda)\cap \dot{\CB}_{\frac{1}{n}}(\dot\mu)\neq\emptyset$ for any $n\in\N$. Picking $\dot\lambda_n\in \dot{\CB}_{\frac{1}{n}}(\dot\lambda)\cap \dot{\CB}_{\frac{1}{n}}(\dot\mu)$ we get a sequence $(\dot\lambda_n)$ that converges to both $\dot\lambda$ and $\dot\mu$. It follows from the characterization of converging sequence given in the previous section that $\int_k d(\breve c(\dot\lambda_n))\longrightarrow\int_k d(\breve c(\dot\lambda))$ for any generic arc $k\subset  S\setminus \lambda^{(p)}$, and $\liminf \int_k d(\breve c(\dot\lambda_n))\geq\pi$ for any generic arc $k$ that intersects $\lambda^{(p)}$ and that the same is true for $\dot\mu$. It follows that $\lambda^{(p)}=\mu^{(p)}$ and $\breve c(\dot\lambda)\cap S\setminus\mu^{(p)}=\breve c(\dot\mu)\cap S\setminus\lambda^{(p)}$, hence $\breve c(\dot\lambda)=\breve c(\dot\mu)$ and $\dot\lambda=\dot\mu$.

Thus $\ML(S)/\CR$ is Hausdorff and since $\overline{\dot{\CB}_{\frac{1}{n+1}}(\dot\lambda)}\subset \dot{\CB}_{\frac{1}{n}}(\dot\lambda)$ for any $\dot\lambda\in\ML(S)/\CR$ and any $n$ large, by \cite[Proposition 1.5.5]{ryszard:general:topology}, $\ML(S)/\CR$ is regular.
\end{proof}

As was already mentioned, $\ML(S)$ is a finite dimensional topological manifold (\cite[Proposition 8.10.5]{thurston:notes}), in particular it is a hereditary Lindel\"{o}f space. By \cite[Theorem 3.8.7]{ryszard:general:topology}) $\ML(S)/\CR$ is also a hereditary Lindel\"{o}f space. Since it is first countable, we deduce the following fact from \cite[Theorems 3.10.1 and 3.10.31]{ryszard:general:topology}):

\begin{lemma}		\label{quotient equivalence}
Given a closed hyperbolic surface $S$, compactness and sequential compactness are equivalent in $\ML(S)/\mathcal R$ (equipped with the tubular topology).
\end{lemma}

This Lemma will allow us to use sequences to prove Theorem \ref{two}.


\subsection{Hyperbolic $3$-manifolds}
Let $M$ be an orientable compact $3$-manifold whose interior admit a complete hyperbolic metric. For most of the paper we will not be interested in the torus component of $\partial M$. For this reason, we will abuse the notation and denote by $\partial M$ the union of the boundary components with negative Euler characteristics (unless stated otherwise).

Let $\sigma$ be a hyperbolic metric on $\IM$. Given an isometry from the interior of $\tilde M$ to $\Hp^3$, the covering transformations yield a discrete faithful representation $\rho:\pi_1(M)\rightarrow Isom^+(\Hp^3)=PSL_2(\C)$. The image  $\rho(\pi_1(M))$ is a finitely generated torsion free Kleinian group and we have an isometry $h_\sigma:(\IM,\sigma)\to\Hp^3/\rho(\pi_1(M))$. The representations which appear in this way will be called {\em representations associated to $\sigma$} and a homeomorphism $(\IM,\sigma)\to\Hp^3/\rho(\pi_1(M))$ isotopic to $h_\sigma$ will be said to be {\em associated to $\sigma$ and $\rho$}. The set of representations associated to $\sigma$ is the set of all representations conjugate to $\rho$. Two metrics that differ by a homeomorphism homotopic to the identity have the same associated representations but they have the same associated homeomorphisms only if they differ by a homeomorphism isotopic to the identity.

The {\em convex core} of $\Hp^3/\rho(\pi_1(M))$ is the quotient by $\rho(\pi_1(M))$ of the smallest convex subset of $\Hp^3$ that is invariant under the action of $\rho(\pi_1(M))$. This set $N(\rho)$ is the image under the isometry $h_{\sigma}$ of the convex core $N(\sigma)$ defined in the introduction. Using the nearest point projection $\IM\to N(\sigma)$ we define an embedding $e_\sigma:N(\sigma)\to M$ that will also be said to be {\em associated to $\sigma$}.


The boundary of $N(\sigma)$ is totally geodesic on the complementary of a geodesic lamination and the amount of bending define a transverse measure on this geodesic lamination. When $\sigma$ is geometrically finite, $e_\sigma(\partial N(\sigma))\subset \partial_{\chi<0} M$ is the complementary in $\partial_{\chi<0} M$ of a family of annuli corresponding to the rank $1$ cusps of $\sigma$. Putting a weight equal to $\pi$ on the core curve of each of these annuli and adding the image under $e_\sigma$ of the measured geodesic lamination quantifying the bending of $\partial N(\sigma)$ we get the {\em bending measured geodesic lamination} of $\sigma$. The bending maps $b:\GF(M)\to \ML(\partial M)$ maps a geometrically finite metric to its bending measured lamination.

\subsubsection{Laminations on the boundary of $3$-manifolds}	\label{lam boundary}

Let us recall the definitions that were given in the introduction.

\begin{definition}	\label{def:pm}
Let $M$ be a compact hyperbolic $3$-manifold and recall that $\partial M$ is the union of the boundary components with negative Euler characteristics. We denote by $\mathcal{P}(M)$ the set of measured geodesic laminations satisfying the following conditions:
\begin{enumerate}[a)]
	\item the weigth of any closed leaf of $\lambda$ is at most $\pi$;
	\item $\exists\eta>0$ such that, for any essential annulus $E$, $i(\partial E,\lambda)\geq\eta$;
	\item $i(\lambda,\partial D)>2\pi$ for any essential disc $D$.
\end{enumerate}
\end{definition}

By \cite{bonahon:otal} and \cite{lecuire:plissage}, $\mathcal{P}(M)$ is the image of the bending map: $\mathcal{P}(M)=b(\GF(M))$. For a compact orientable hyperbolic $3$-manifold $M$, we have $\breve c(\mathcal P(\partial M)/\CR)=\mathcal P(\partial M)$.

In the introduction we also mentioned the set $\mathcal{D}(M)$ of doubly incompressible laminations:

\begin{definition}
A measured geodesic lamination $\lambda\in\mathcal{ML}(\partial M)$ is {\em doubly incompressible} if  $\exists\eta>0$ such that $i(\partial E,\lambda)>\eta$ for any essential annulus or disc  $E$.
\end{definition}

The inclusion $\mathcal{P}(M)\subset \mathcal{D}(M)$ follows from the definitions. It is also straightforward to see that if $\lambda\in \mathcal{D}(M)$ then $\mu\in\mathcal{D}(M)$ for any measured lamination $\mu$ with $\mu\mathcal{R}\lambda$.

\subsubsection{Limits of diverging sequences of representations}

In this section we will give the essential definitions and results from Thurston's compactification of Teichm\"uller space and Morgan-Shalen compactification of character varieties that will be used in this article.

From Thurston's compactification of Teichm\"uller space we will only use the the last part of \cite[theorem 2.2]{thurston:hypII}.

\begin{lemma}	\label{thurston}
Given a finite type surface $S$, there is a constant $C$ with the following property: for any complete hyperbolic metric $s$ on $S$, there is a measured geodesic lamination $\mu(S)$ such that
	\begin{equation}
		i(\mu(S),c)\leq\ell_{s}(c)\leq i(\mu(S),c)+ C\ell_{s_0}(c).	 \label{thucom}
	\end{equation}
for any simple close curve $c\in S$.
\end{lemma}

The metric $s_0$ is an arbitrarily chosen reference metric.
The length $\ell_{s}(c)$ is the length of the unique $s$-geodesic in the free homotopy class of $c$. This length function naturally extends to weighted multi-curves and then to measured geodesic laminations by taking limits. Then equality (\ref{thucom}) still holds when $c$ is any measured geodesic lamination.

A proof of this Lemma can be found in \cite[Exposé 8, \textsection II.1]{FLP} (combined with \cite{levitt:foliations:laminations}).\\

We will also use Culler-Morgan-Shalen's compactification of character varieties by actions on real trees. Let $M$ be a compact $3$-manifold and $\rho_n:\pi_1(M)\to PSL(2,\C)$ a sequence of discrete and faithful representations. We say that $(\rho_n)$ tends to a minimal action of $\pi_1(M)$ on an $\R$-tree $\mathcal{T}$ when there exists a sequence $(\varepsilon_n)$ of positive numbers tending to $0$ and a minimal action of $\pi_1(M)$ on an $\R$-tree $\mathcal{T}$ such that, for any $g\in\pi_1(M)$, we have $\ell_\mathcal{T}(g)=\lim_{n\longrightarrow\infty} \varepsilon_n\ell_{\rho_n}(g)$. Here $\ell_{\rho_n}(g)$ is the distance of translation of $\rho_n(g)$ and $\ell_\mathcal{T}(g)$ is the distance of translation of $g$ acting on $\mathcal{T}$.
The theory of \cite{morgan:shalen:degenerationsI} establishes the following:

\begin{theorem}
Let $(\rho_n)$ be a sequence of faithful discrete representations of $\pi_1(M)$ for a compact $3$-manifold $M$. Either there is a subsequence indexed by $k:\N\to\N$ and $q_n\in PSL(2,\C)$ such that $q_{k(n)}\rho_{k(n)}(g)q_{k(n)}^{-1}$ converges for any $g\in\pi_1(M)$ or there is a subsequence indexed by $k:\N\to\N$ such that $(\rho_{k(n)})$ tends to a minimal action of $\pi_1(M)$ on a non-trivial $\R$-tree.
\end{theorem}

\section{Convergence of representations}

As a first step in the proof, we will show that given a sequence of metrics with bending laminations converging modulo $\mathcal{R}$, there is a sequence of associated representations containing a converging subsequence. 

\begin{prop}		\label{algebrique}
	Let $M$ be a compact, orientable, hyperbolic 3-manifold with boundary and let $(\sigma_n)$ be a sequence of geometrically finite metrics on the interior of $M$. Let $\lambda_n\subset\mathcal{ML}(\partial M)$ be the bending measured geodesic lamination of $\sigma_n$. Suppose that $\dot\lambda_n\in\mathcal{ML}(\partial M)/\mathcal{R}$ tends to  $\dot\lambda_\infty\subset \mathcal{P}(M)/\mathcal{R}$. Then, up to extracting a subsequence, there is a sequence of representations $(\rho_n)$  such that $\rho_n$ is associated to $\sigma_n$ for any $n$ and that $(\rho_n(g))$ converges for any $g\in\pi_1(M)$.
\end{prop}

The starting point of the proof is an upper bound on the length of the bending lamination that was first observed by Bonahon-Otal (\cite[Lemma 9]{bonahon:otal}). Our case is slightly more general but can be deduced from  \cite[Lemma 4.1]{lecuire:plissage}.

\begin{lemma}	\label{borne}
	Under the hypothesis of Proposition \ref{algebrique} the lengths $\ell_{\sigma_n}(\lambda_n)$ are uniformly bounded.
\end{lemma}

\begin{proof}
	The proof of \cite[Lemma 4.1]{lecuire:plissage} extends straightforwardly to this situation.	
\end{proof}

Notice that an explicit bound is provided in \cite{bridgeman:canary:bounding:bending}.

From this point, the proof of the convergence of representations follows a fairly extensive line of works initiated with Thurston's Double Limit Theorem (\cite[Theorem 4.1]{thurston:hypI}) and continued in its diverse generalizations, in particular by Canary \cite{canary:schottky}, Otal \cite{otal:schottky}, Kleineidam-Souto \cite{kleineidam:souto:function}, Lecuire \cite{lecuire:masur} and Kim-Lecuire-Ohshika \cite{klo:convergence}. It is slightly disappointing that none of these generalizations directly applies to our case where some leaves might have an infinite weight. In order to still make the best use of those precursors, we will start with the following result which stems from Otal's work (\cite{otal:fibre}) and is essentially a combination of \cite[Proposition 6.1]{lecuire:masur} and \cite[Theorem 6.5]{lecuire:masur}.

\begin{theorem}	\label{contin}
Let $M$ be an orientable compact hyperbolic $3$-manifold and let $(\rho_n)$ be a sequence of discrete and faithful representations of $\pi_1(M)$ into $PSL_2(\C)$ tending to a minimal action of $\pi_1(M)$ on a non-trivial $\R$-tree $\mathcal{T}$. Let $\lambda$ be a measured geodesic lamination in $\mathcal{D}(M)$ and let $E$ be a geodesic lamination containing $|\lambda|$, such that $S(E)=S(|\lambda|)$. Consider $\varepsilon_n\longrightarrow 0$ such that $\forall g\in\pi_1(M)$, we have $\varepsilon_n\ell_{\rho_n}(g)\longrightarrow\ell_\mathcal{T}(g)$. Then there is at least one connected component $E_1$ of $E$ such that:
	\begin{enumerate}[1)]
		\item either $E_1$ contains an irrational geodesic lamination;
		\item or $E_1$ is a closed leaf and $E_1$ does not lie in $\partial S(E_i)$ for any component $E_i\neq E_1$ of $E$.
	\end{enumerate}
and there exists a neighborhood $\mathcal{N}(E_1)$ of $E_1$ in the space of geodesic laminations and constants $K,n_0$ such that for any simple closed curve $c\subset \mathcal{N}(E_1)$ and for any $n\geq n_0$,
	$$\varepsilon_n\ell_{\rho_n}(c)\geq K \ell_{s_0}(c).$$
\end{theorem}

Again, $s_0$ is an arbitrary chosen reference metric.

\begin{proof}
By combining \cite[Proposition 6.1]{lecuire:masur} (generalizing \cite[Theorem 3.1.4]{otal:fibre}) and \cite[Theorem 6.5]{lecuire:masur}(see also \cite[Theorem 3.7]{lecuire:double} generalizing \cite[Theorem 4.0.1]{otal:fibre}) we get that there is a component $E_1$ of $|\lambda|$ that satisfies the last part. To upgrade $E_1$ to a component of $E$ of type $1)$ or $2)$ we need to have a more careful look at the proof of \cite[Proposition 6.1]{lecuire:masur}. Since $S(E)=S(|\lambda|)$, there is a sequence of weighted multi-curves $\lambda_n$ that converges in $\mathcal{ML}(\partial M)$ to a measured lamination supported by $\lambda$ and whose supports converge to $E$ in the Hausdorff topology. Building on the works of Morgan and Otal, \cite{morgan:otal} and \cite{otal:fibre}, we construct for each $n$ a measured geodesic lamination $\beta_n$ together with, for each component $S$ of $\partial M$, a map $\tilde S\to \mathcal T$ which is injective on the preimage of every leaf of $\lambda_n$ that intersects $|\beta_n|$ transversely. Then we extract a subsequence such that $(|\beta_n|)$ converges in the Hausdorff topology to a lamination $B$. By \cite[Lemma 6.2]{lecuire:masur}, $E$ intersects $B$ transversely. It follows that there is a component $E_1$ of $E$ that satisfies $1)$ or $2)$ and intersects $B$ transversely. The arguments used to conclude the proof of \cite[Proposition 6.1]{lecuire:masur} show that $E_1$ is realized in $\tau$ and the conclusion follows from  \cite[Theorem 6.5]{lecuire:masur}.
\end{proof}

From there we are going to do cut and paste operations on $\lambda_n$ to obtain the geodesic lamination $E$ from Theorem \ref{contin} as a limit of weighted multi-curve with controlled lengths. This control will fend off the situation described in the conclusion of Theorem \ref{contin} and hence prevent any subsequence of $(\rho_n)$ to degenerate according to Culler-Morgan-Shalen theory.


\subsection{Cut and paste operations}

We will now build the components of $E$ corresponding to case $1)$ of Theorem \ref{contin}.

\begin{lemma}	\label{elle}
	Let $S$ be a finite type surface and $s_n$ a sequence of complete hyperbolic metric on $S$. Consider a sequence $(\lambda_n)\subset \ML(S)$ such that $\dot\lambda_n\in\mathcal{ML}(S)/\mathcal{R}$ tends to $\dot\lambda$. Let $\gamma$ be a minimal sublamination of $\lambda$ that is not a supported on a closed curve. Then there is a sequence of weighted multi-curves $\gamma_n$ such that:
	\begin{enumerate}[ - ]
		\item $|\gamma_n|\subset S(\gamma)$,
		\item $\ell_{s_n}(\gamma_n)\leq\ell_{s_n}(\lambda_n)+i(\lambda_n,\partial S(\gamma))\ell_{s_n}(\partial S(\gamma))$,
		\item a subsequence of $\{\gamma_n\}$ converges to a non-trivial measured geodesic lamination with the same support as $\gamma$.
	\end{enumerate}
\end{lemma}

\begin{proof}
	Since weighted simple closed curves are dense in $\ML(S)$ (\cite[Theorem 3.1.3]{penner:harer}), we can approximate $\lambda_n$ by weighted simple closed curves $\nu_n$ so that $\dot\nu_n$ converges to $\dot\lambda$ and that $i(\nu_n,\lambda_n)\longrightarrow 0$. Up to slightly reducing the weight $p_n$ of $\nu_n$, we may assume that $\ell_{s_n}(\nu_n)\leq\ell_{s_n}(\lambda_n)$ and $i(\nu_n,\partial S(\gamma))\leq i(\lambda_n,\partial S(\gamma))$ for any $n$. If $|\nu_n|\cap S(\lambda)$ is a simple closed curve, we are done, taking $\gamma_n=\nu_n$. Otherwise $\gamma_n$ will be a component of the boundary of a neighborhood of $k_n\cup \partial S(\gamma)$ for a carefully chosen component $k_n$ of $|\nu_n|\cap S(\gamma)$.
	
	We are going to progressively be more selective on the choice of $k_n$ to ensure that $\gamma_n$ has the expected properties. Before that, we set up some notations. For a component $k_n$ of $|\nu_n|\cap S(\gamma)$, we denote by $\eps(k_n)\in\{1,2\}$ the number of components of $\partial S(\gamma)$ containing the endpoints of $k_n$ and by $\CN(k_n)$ a neighborhood of the union of $k_n$ and those $\eps(k_n)$ components. We denote by $[k_n]$ the union of the components of $|\nu_n|\cap S(\gamma)$ in the isotopy class of $k_n$ relative to $\partial S(\gamma)$. Notice that even though there might be a growing number of components of $|\nu_n|\cap S(\gamma)$, there is an upper bound, only dependent on $S(\gamma)$, on the number of isotopy classes. Let $v(k_n)$ be the number of components of $[k_n]$.
	
	Given $k_n$, since $S(\gamma)$ is not a pair of pants, there is at least one component $g_n$ of $\partial\CN(k_n)\cap S(\gamma)$ that is not peripheral. Equipping $g_n$ with a weight equal to $\frac{p_nv(k_n)}{\eps(k_n)}$ (remember that $p_n$ is the weight of $\nu_n$) we get a measured geodesic lamination $\gamma_n$.
	
	By construction, we have $\ell_{s_n}(\gamma_n)\leq p_n v(k_n)(\ell_{s_n}(k_n)+\ell_{s_n}(\partial S(\gamma)))$. If we choose for $k_n$ the arc in $[k_n]$ with the smallest $s_n$ length, then $p_n v(k_n)\ell_{s_n}(k_n)\leq p_n\ell_{s_n}([k_n])\leq\ell_{s_n}(\nu_n)\leq\ell_{s_n}(\lambda_n)$. On the other hand $p_n v(k_n)\leq p_n\sharp\{|\nu_n|\cap \partial S(\gamma)\}=i(\nu_n,\partial S(\gamma))\leq i(\lambda_n,\partial S(\gamma))$. Putting these inequalities together, we get $\ell_{s_n}(\gamma_n)\leq\ell_{s_n}(\lambda_n)+i(\lambda_n,\partial S(\gamma))\ell_{s_n}(\partial S(\gamma))$.
	
	To ensure the last property, pick a simple closed curve $c\in S(\gamma)$ that is not isotopic to a boundary component and for each $n$ consider an isotopy class $[k_n]$ that maximizes the cardinality $\sharp\{[k_n]\cap c\}$. Then $\frac{\sharp\{[k_n]\cap c\}}{\sharp\{|\nu_n|\cap c\}}$ is bounded away from $0$. Let $k_n$ be the arc in $[k_n]$ with the smallest $s_n$ length. Since $c$ is not isotopic to a boundary component, there is a component $g_n$ of $\partial\overline{\mathcal{N}(k_n)}$ that satisfies $i(g_n,c)=\eps(k_n) \sharp\{k_n\cap c\}$. When we equip $g_n$ with a weight equal to $\frac{p_nv(k_n)}{\eps(k_n)}$ to produce $\gamma_n$, we then get $i(\gamma_n,c)=p_n\sharp\{[k_n]\cap c\}$. Since $\frac{\sharp\{[k_n]\cap c\}}{\sharp\{|\nu_n|\cap c\}}$ is bounded away from $0$ and $p_n\sharp\{|\nu_n|\cap c\}=i(\nu_n,c)$, then $i(\gamma_n,c)$ is bounded away from $0$ and bounded. In particular, a subsequence of $\{\gamma_n\}$ converges to a non trivial measured geodesic lamination $\gamma_\infty$. By construction $i(\gamma_n,\nu_n)\leq p_nv(k_n) i(\nu_n,\partial S(\gamma))\leq i(\nu_n,\partial S(\gamma))^2\longrightarrow 0$. Since $\dot\nu_n$ converges to $\dot\lambda_\infty$, it follows that $\gamma_\infty\in\ML(S(\gamma))$ has the same support as $\gamma$.
\end{proof}

\subsection{Comparing the length of closed curves}

Next we explain how to control the lengths of the closed leaves of $\dot\lambda_\infty$ corresponding to the case $2)$ of Theorem \ref{contin}.

\begin{lemma}	\label{ce}
	Let $M$ be a compact, orientable, hyperbolic 3-manifold with boundary and let $(\sigma_n)$ be a sequence of geometrically finite metrics on $\mathring{M}$. We denote by $\lambda_n\subset\mathcal{ML}(\partial M)$ the bending measured geodesic lamination of $\sigma_n$ and by $s_n$ be the hyperbolic  metric on $\partial M$ induced by an embedding $e_n:N(\rho_n)\to M$ associated to $\sigma_n$.
	Suppose that $\dot\lambda_n$ converges to $\dot\lambda_\infty$. Let $d$ be a closed leaf of $\breve c(\dot\lambda_\infty)$ that does not lie in the boundary of the minimal supporting surface $\partial S(E_j)$ of any component $E_j\neq d$ of $\breve c(\dot\lambda_\infty)$. If $\ell_{s_n}(\lambda_n)$ is bounded but $\ell_{s_n}(d)\longrightarrow\infty$ then there is an essential subsurface $F\subset \partial M$ such that:
	\begin{enumerate}[- ]
		\item $|\dot\lambda_\infty|\cap\mathrm{int}(F)=d$,
		\item $d$ is not isotopic in $\partial M$ to a component of $\partial F$,
		\item for any simple closed curve $c\subset F$ that intersects $d$ transversely, we have $$\frac{\ell_{\sigma_n}(d)}{\ell_{\sigma_n}(c)}\longrightarrow 0$$
	\end{enumerate}
\end{lemma}

Recall from \textsection \ref{lam boundary} that $\breve c(\dot\lambda_n)$ is the representative of $\dot\lambda_n$ without any leaf with a weight greater than $\pi$.

\begin{proof}
	We let $F$ be the largest closed subsurface of $\partial M$, up to isotopy, containing $d$ such that $i(c,\breve c(\dot\lambda_\infty))=0$ for any simple closed curve $c\subset F\setminus  d$. We have $\breve c(\dot\lambda_\infty)\cap\mathrm{int}(F)=d$, and, since $d$ does not lie in the boundary of the minimal supporting surface of any other component of $\breve c(\dot\lambda_\infty)$, $d$ is not isotopic to a component of $\partial F$.
	
	Let $c\subset F$ be a simple closed curve that intersects $d$ transversely. We are first going to show that $\frac{\ell_{s_n}(d)}{\ell_{s_n}(c)}\longrightarrow 0$. Then we will show that $c$ is only slightly bend to obtain a constant $K$ such that $\ell_{s_n}(c)\leq K\ell_{\sigma_n}(c)$.
	
	\begin{claim}	\label{compare}
		$\frac{\ell_{s_n}(d)}{\ell_{s_n}(c)}\longrightarrow 0$
	\end{claim}
	
	\begin{proof}
		We are going to show that any subsequence contains a further subsequence for which the conclusion is satisfied. Let $\mu(s_n)$ be a sequence of measured geodesic laminations associated to $s_n$, i.e. such that $i(\mu,\mu(s_n))\leq\ell_{s_n}(\mu)\leq i(\mu,\mu(s_n))+ C\ell_{s_0}(\mu)$ for any $\mu\in\mathcal{ML}(\partial M)$ as in Lemma \ref{thurston}. Since weighted multi-curves are dense in $\mathcal{ML}(\partial M)$, we may assume that $\mu(s_n)$ is a weighted multi-curve for any $n$. Fix a hyperbolic metric $s_0$ on $\partial M$, consider an arc $\kappa\subset c$ with one endpoint in $d$ which is small enough so that any simple $s_0$-geodesic intersecting $\kappa$ intersects $d$ and set $K_n=\frac{\sharp \{|\mu(s_n)|\cap c\}}{\sharp\{|\mu(s_n)|\cap d\}}$. Extract a subsequence such that $(|\mu(s_n)|)$ converge in the Hausdorff topology to a geodesic laminations $m_\infty$.  Since $\ell_{s_n}(d)\longrightarrow\infty$, $i(\mu(s_n),d)\longrightarrow\infty$. On the other hand $d$ is a closed leaf of $\dot\lambda_\infty$ and $i(\mu(s_n),\lambda_n)$ is bounded. This is possible only if $m_\infty$ contains $d$ and leaves spiraling around $d$ hence $\mu(s_n)$ spirals more and more around $d$. It follows that we have $K_n\longrightarrow\infty$ and $i(\mu(s_n),c)\geq K_n i(\mu(s_n),d)$. This leads to $\ell_{s_n}(c)\geq K_n(\ell_{s_n}(d)-C\ell_{s_0}(d))$ which concludes the proof since $\ell_{s_n}(d)\longrightarrow\infty$ by assumption.
	\end{proof}

	To conclude the proof of Lemma \ref{ce}, we need to compare $\ell_{s_n}(c)$ and $\ell_{\sigma_n}(c)$. This will follow from the fact that $c$ can be cut into arcs with small bending and unbounded length.
	\begin{claim}	\label{sbend}
		For any sequence of arcs $\kappa_n\subset c$ such that $\partial\kappa_n\subset\lambda_n$ and $\underline\lim\; \int_{\kappa_n} d\lambda_n>0$, we have $\ell_{s_n}(\kappa_n)\longrightarrow \infty$.
	\end{claim}
	\begin{proof}
		We are going to compare the behavior of the sequences $(s_n)$ and $(\lambda_n)$ near $d$. If $d$ is a leaf of $\lambda_n$, then its weight in $\lambda_n$ converges to its weight in $\breve c(\dot\lambda_\infty)$. Since $\ell_{s_n}(d)$ tends to $\infty$ and $\ell_{s_n}(\lambda_n)$ is bounded, we have $d\not\subset|\lambda_n|$ for large $n$.
		
		Let us consider a subsequence such that $(|\lambda_n|)$ converges to some geodesic lamination $ L_\infty$ in the Hausdorff topology. Since $d$ is a leaf of $c(\lambda_\infty)$, $d\subset L_\infty$ and $i(\lambda_n,d)\longrightarrow 0$. Let $\mathcal{N}(d)$ be an annular neighborhood of $d$. Since $d$ is not a leaf of $\lambda_n$ for large $n$, then $\mathcal{N}(d)\cap|\lambda_n|$ is the union of disjoint segments joining the two components of $\partial\mathcal{N}(d)$ and turning many times arouns $d$ (see figure \ref{voisin}). In particular any leaf of $\lambda_n$ entering $\mathcal{N}(d)$ intersects $d$.
		\begin{figure}[hbtp]
			\psfrag{a}{$c$}
			\psfrag{b}{$d$}
			\psfrag{c}{$\lambda_n$}
			\psfrag{d}{$x_n$}
			\psfrag{e}{$y_n$}
			\psfrag{f}{$\kappa_n$}
			\psfrag{g}{$c_{x_n y_n}$}
			\centerline{\includegraphics{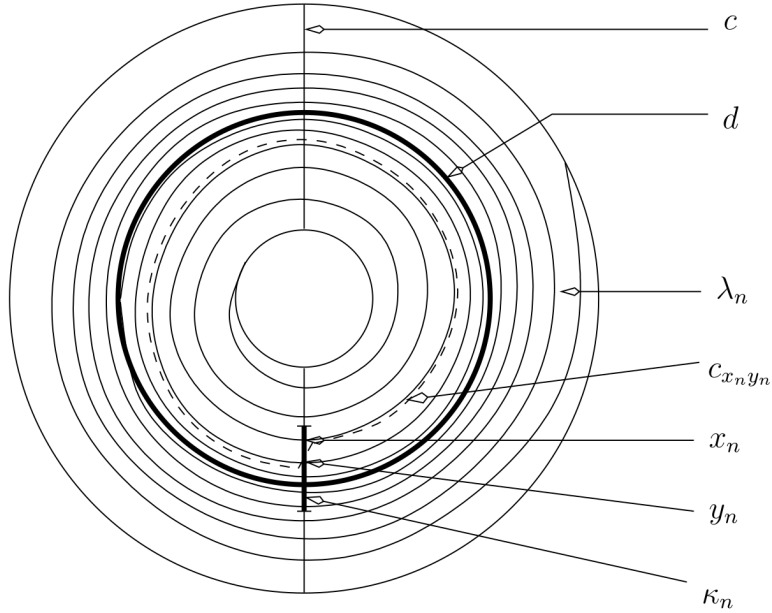}}
			\caption{Picture of $\mathcal{N}(d)$}
			\label{voisin}
		\end{figure}

		Let $\kappa_n\subset c$ be a sequence of arcs such that $\partial\kappa_n\subset\lambda_n$ and $\underline\lim\; \int_{\kappa_n} d\lambda_n>0$. Since $F\cap |\dot\lambda_\infty|=d$,  $\int_{c\setminus \mathcal{N}(d)} d\lambda_n$ tends to $0$. Therefore it is sufficient to prove Claim \ref{sbend} for $\kappa_n\subset \mathcal{N}(d)$. Given an orientation of $c$  we give an order induced from the orientation on the
		points of $\kappa_n\cap\lambda_n$. Since $i(\lambda_n,d)\longrightarrow 0$, for $n$ large enough, $i(\lambda_n,d)<\frac{\int_{\kappa_n} d\lambda_n}{2}$.
		Since any leaf of $\lambda_n$ entering $\mathcal{N}(d)$ intersects $d$, for large $n$, $\kappa_n$ contains an arc $\kappa'_n$ such that $\int_{\kappa'_n} d\lambda_n\geq\frac{ \int_{\kappa_n} d\lambda_n}{2}$ and such that any point $x_n$ of $\kappa'_n\cap\lambda_n$ is joined to a point $y_n\neq x_n$ of $\kappa_n \cap\lambda_n$ by an arc $]x_n,y_n[\subset\lambda_n\setminus \kappa_n$. Fix $n$ and consider a point $x_n$ in $\kappa'_n\cap\lambda_n$ and the associated point $y_n$. Let $c_{x_n y_n}$ denote the simple closed curve which is the union of the arc $[x_n, y_n]\subset\lambda_n$ and of the arc $\kappa_{x_n y_n}\subset\kappa_n$ joining $x_n$ to $y_n$ (see figure \ref{voisin}). This curve $c_{x_n y_n}$ is homotopic to $d$, hence $\ell_{s_n}(c_{x_n y_n})$ is greater than $\ell_{s_n}(d)$. It follows that $\forall n$, $\forall x_n\subset\kappa'_n\cap\lambda_n$, we have $\ell_{s_n}([x_n,y_n])\geq\ell_{s_n}(d)-\ell_{s_n}(\kappa_n)$. Since $\ell_{s_n}(\lambda_n)\geq (\int_{\kappa'_n} d\lambda_n) \inf\{\ell_{s_n}([x_n,y_n])\,|\,x_n\subset\kappa'_n\cap\lambda_n \}$, we have $\ell_{s_n}(\lambda_n)\geq (\int_{\kappa'_n} d\lambda_n )(\ell_{s_n}(d)-\ell_{s_n}(\kappa_n))$. Thus, we get $\ell_{s_n}(\kappa_n)\geq\ell_{s_n}(d)-\frac{2}{\int_{\kappa_n} d\lambda_n}\ell_{s_n}(\lambda_n)$. Since $\ell_{s_n}(d)\longrightarrow\infty$ and $\ell_{s_n}(\lambda_n)$ is bounded, this leads to  $\ell_{s_n}(\kappa_n)\longrightarrow\infty$.
	\end{proof}
	
	\indent Let $0<\eps<\frac{i(c,c(\lambda))}{2}$. It follows from the description of $\lambda_n\cap\mathcal{N}(E_i)$ that, for large $n$, we may divide $c$ into $\lfloor\frac{i(c,\lambda_n)}{\eps}\rfloor$ segments $\kappa_n^j$ such that $\forall j\leq\lfloor\frac{i(c,\lambda_n)}{\eps}\rfloor$, we have $\eps\leq \int_{\kappa_n^j} d\lambda_n< 2\eps$. By Claim \ref{sbend}, $\forall j$, $\ell_{s_n}(\kappa_n^j)\longrightarrow\infty$.  Given $\eps$ such that $2\eps\leq\frac{\pi}{3}$, let $c_n$ be the union of the geodesic segments of $(M,\sigma_n)$ joining the endpoints of the segments $\kappa_n^j$. By (\cite[Lemma A.2]{lecuire:plissage}), $\exists C_{\eps}$ such that, $\ell_{s_n}(c)\leq C_{\eps}\ell_{\sigma_n}(c_n)$. Moreover, by the Gauss-Bonnet formula, the geodesic arcs comprising $c_n$ have interior angles greater than $\pi-2\frac{\pi}{3}=\frac{\pi}{3}$. Hence the curve $c_n$ is the union of long segments with incidence angles greater than $\frac{\pi}{3}$. In this situation, it is a classical result (compare with the proof of \cite[Lemma A.1]{lecuire:plissage}) that $\exists K_{\eps}$ such that we have $\ell_{\sigma_n}(c_n)\geq K_{\eps}\ell_{\sigma_n}(c)$. 
	Now we get $\ell_{s_n}(c)\leq K_{\eps} C_{\eps}\ell_{\sigma_n}(c)$. Combined with Claim \ref{compare}, this inequality concludes the proof of Lemma \ref{ce}.
\end{proof}

\subsection{Proof of the convergence of representations}
We are now ready to prove Proposition \ref{algebrique}.

\begin{proof}[Proof of Proposition \ref{algebrique}]
We say that a sequence of representations $(\rho_n)$ converges if $(\rho_n(g))$ converges in $PSL_2(\C)$ for any $g\in\pi_1(M)$. Assume that the conclusion of Proposition \ref{algebrique} does not hold, i.e. no sequence of associated representations $(\rho_n)$ contains a converging subsequence. 
By Culler-Morgan-Shalen's theory, up to extracting a subsequence, there is $\eps_n\longrightarrow 0$ and a minimal action of $\pi_1(M)$ on a non-trivial $\R$-tree $\mathcal{T}$ such that $\varepsilon_n\ell_{\rho_n}(g)\longrightarrow\ell_\mathcal{T}(g)$ for any $g\in\pi_1(M)$. 

For each minimal sublamination $\gamma_i$ of $\dot\lambda_\infty$ that is not supported on a closed curve, we have constructed in Lemma \ref{elle} a sequence of weighted multi-curves $\gamma_n$ with controlled lengths. Extract a subsequence such that $|\gamma_n|$ converges in the Hausdorff topology to a lamination $E_i$. By Lemma \ref{elle}, $E_i$ contains $|\gamma_i|$ and $E_i\subset S(\gamma_i)$ but we have not ruled out the possibility that $E_i$ has some non-recurrent leaves that do not belong to $|\gamma_i|$. Adding to all those Hausdorff limits $E_i$ the simple closed curves in $|\dot\lambda_\infty|$, we get a geodesic lamination $E$ such that $|\dot\lambda_\infty|\subset E$ and $S(E)=S(\dot\lambda_\infty)$. All the hypothesis of Theorem \ref{contin} are satisfied, thus we have a connected component $E_1$ of $E$, a neighborhood $\mathcal{N}(E_1)$ of $E_1$ in the space of geodesic laminations and constants $K,n_0$ such that for any simple closed curve $c\subset \mathcal{N}(E_1)$ and for any $n\geq n_0$,
$\varepsilon_n\ell_{\sigma_n}(c)\geq K \ell_{s_0}(c)$.

If $E_1$ is a not a simple closed curve then by Lemma \ref{elle}, we have $\ell_{s_n}(\gamma_n)\leq\ell_{s_n}(\lambda_n)+i(\lambda_n,\partial S(\gamma))\ell_{s_n}(\partial S(E_1))$. Since $S(E)=S(|\dot\lambda_\infty|)$, $i(\lambda_n,\partial S(E_1))\longrightarrow 0$. It follows from \cite[Lemma A.1]{lecuire:plissage} that there are $A_n\longrightarrow 0$ and $C_n\longrightarrow 1$ such that $\ell_{s_n}(\partial S(E_1))\leq C_n(\ell_{\rho_n}(\partial S(E_1))+A_n)$. Since $\eps_n\ell_{\rho_n}(\partial S(E_1))\longrightarrow \ell_{\mathcal{T}}(\partial S(E_1))$, $i(\lambda_n,\partial S(E_1))\longrightarrow 0$ and $(\ell_{s_n}(\lambda_n))$ is bounded (Lemma \ref{borne}), $\eps_n\ell_{\rho_n}(\gamma_n)\leq\eps_n\ell_{s_n}(\gamma_n)\longrightarrow 0$, contradicting the conclusion of Theorem \ref{contin}.

If $E_1$ is a closed curve, then by $2)$ of Theorem \ref{contin}, it does not lie on the boundary of a the minimal supporting surface of another component of $E$. By Lemma \ref{ce}, there is a simple closed curve $c$ such that $\frac{\ell_{\rho_n}(E_1)}{\ell_{\rho_n}(c)}\longrightarrow 0$. On the other hand, we have $\eps_n\ell_{\rho_n}(c)\longrightarrow \ell_\mathcal{T} (c)$. Thus we get $\eps_n\ell_{\rho_n}(E_1)\longrightarrow 0$, contradicting again the conclusion of Theorem \ref{contin}.

We have shown that, assuming that the conclusion of Proposition \ref{algebrique} does not hold leads to a contradiction between Theorem \ref{contin} and Lemma \ref{elle} or \ref{ce}. This concludes the proof.
\end{proof}

\section{Upgrading the convergence}

 Let $M$ be a compact orientable hyperbolic 3-manifold, let $(\sigma_n)$ be a sequence of geometrically finite metrics on $\mathring{M}$ and let $\lambda_n\subset\mathcal{ML}(\partial M)$ be the bending measured geodesic lamination of $\sigma_n$. Assume that $\dot\lambda_n$ converges to $\dot\lambda_\infty\in\mathcal{P}( M)/\mathcal{R}$, by Lemma \ref{algebrique}, there is a sequence $\rho_n:\pi_1(M)\rightarrow Isom(\Hp^3)$ of representations associated to $(\sigma_n)$ such that a subsequence of $(\rho_n)$ tends to a representation $\rho_{\infty}:\pi_1(M)\rightarrow Isom(\Hp^3)$. To simplify the notations, assume that the whole sequence $(\rho_n)$ converges to $\rho_\infty$. By \cite{jorgensen:mobius} (see also \cite{chuckrow}), $\rho_{\infty}$ is discrete and faithful. In this section, we are going to show that $\rho_{\infty}$ is geometrically finite and associated to a metric on the interior of $M$.
 
 We start with a variation on the classical fact that weighted multi-curves are dense in $\ML(S)$. 

 \begin{lemma}		\label{pantssmallmeasure}
 	Let $S$ be a closed surface and $\lambda\in\ML (S)$. For any $\eps>0$, there is a pants decomposition $P$ of $S$ such that $i(\lambda,P)\leq\eps$. Furthermore, for any minimal sublamination $\mu$ of $\lambda$, either the support of $\mu$ is a simple closed curve or there is a simple closed curve $c_\mu$ such that $i(\lambda, c_\mu)\leq\eps$ and $i(P,d)+i(c_\mu,d)>0$ for any non-peripheral simple closed curve $d\subset S(\mu)$.
 \end{lemma}
 
 We will call $c_\mu$ a {\em transverse} for $P\cap S(\mu)$.
 
 \begin{proof}
 Start with a multicurve $P$ made up of all the closed leaves of $\lambda$. If there is a simple closed curve $c\subset S\setminus  P$ not isotopic to a leaf of $P$ and satisfying $i(c,\lambda)=0$, then we add $p$ to $P$. We repeat this process until any curve $c\subset S\setminus  P$ with $i(c,\lambda)=0$ is isotopic to a leaf of $P$. Consider $\eta>0$ that will be fixed later and pick a component $F$ of $S\setminus  P$. It follows from the construction that if $F$ is not a pair of pants then $|\lambda|\cap F$ is minimal and filling. We will recursively build a pants decomposition of $F$. In the first step, we construct, following \cite[I.4.2.15.]{canary:epstein:green}, a simple closed curve $c_1\subset F$  such that $i(c_1,\lambda)\leq\eta$. Then assume that we have a multicurve $C_j\subset F$ with $j$ leaves such that $i(c,\lambda)\leq 2^{j-1}\eta$ for any leaf $c$ of $C_j$. Consider a component $F_j$ of $F\setminus  C_j$ that is not a pair of pants, let $\kappa$ be a component of $|\lambda|\cap F_j$ and denote by $\mathcal{N}_j\subset F$ a neighborhood of $\kappa\cup(\bar F_j\setminus  F_j)$. At least one component $c_{j+1}$ of $\partial\overline{\mathcal{N}}_j\cap F_j$ is a closed curve that is not isotopic to a component of $\bar F_j\setminus  F_j$. By construction, we have $i(c_{j+1},\lambda)\leq 2^j\eta$. We repeat this construction until we obtain a pants decomposition of $F$. Proceeding in the same way on any component of $S\setminus  P$ that is not a pair of pants, we build a pants decomposition $P$ such that $i(\lambda,c)\leq 2^K\eta$ for any leaf $c$ of $P$ where $K=-\frac{3}{2}\chi(S)$ is the number of curves in a pants decomposition of $S$. We simply need to choose $\eta$ such that $K2^K\eta\leq\eps$ to conclude the construction of $P$.
 	
 	Let $\gamma$ be a minimal sublamination of $\lambda$ and assume that the support of $\gamma$ is not a simple closed curve. Following \cite[I.4.2.15.]{canary:epstein:green} again, consider a sequence of non-peripheral simple closed curves $c_k\subset S(\gamma)$ such that $i(c_k,\gamma)\longrightarrow 0$. For $k$ large enough  $i(c_k,\gamma)\leq\eps$ and $i(P,d)+i(c_k,d)>0$ for any non-peripheral simple closed curve $d\subset S(\gamma)$. The first property is obvious, to show the second property let us notice that any subsequence of $\{c_k\}$ contains a further subsequence that converges in the Hausdorff topology to a lamination $C$. Since $i(c_k,\gamma)\longrightarrow 0$, $C$ contains the support of $\gamma$ and hence intersects any non-peripheral simple closed curve on $S(\gamma)$. It follows that $c_k$ intersects each component of $S(\gamma)\cap P$ for $k$ large enough. In particular, if we take $c_\mu=c_k$ for a large enough $k$, and $d\subset S(\gamma)$ is a non-peripheral simple closed curve such that $i(d,P)=0$ then $d$ is a leaf of $P$ and $i(d,c_\mu)\geq 1$.
 \end{proof}

\subsection{The metrics on the boundary}

We will use Lemma \ref{pantssmallmeasure} and \cite[Lemme A.1]{lecuire:plissage} to compare the induced metric on the boundary of the convex core with the metric on its interior. This will allow us to deduce the convergence of the induced metrics from the convergence of the associated representations. We consider the assumptions that were made at the beginning of \textsection 4, namely $(\dot\lambda_n)$ converges to $\dot\lambda_\infty\in\mathcal{P}( M)/\mathcal{R}$ and $(\rho_n)$ converges to $\rho_\infty$. Let us denote by $\lambda_n^{(p)}$, resp. $\lambda_\infty^{(p)}$, the union of the leaves of $\lambda_n$, resp. $\breve c\dot(\lambda_\infty)$, with a weight equal to $\pi$. By Claim \ref{larger than pi}, for $n$ large enough, $\lambda_n^{(p)}\subset \lambda_\infty^{(p)}$. It is now easy to extract a subsequence such that $\lambda_n^{(p)}$ does not depend on $n$ and set $\lambda^{(p)}=\lambda_n^{(p)}$. Let $e_n: N(\rho_n)\to M$ be an embedding associated to $\sigma_n$ and $\rho_n$. The image of $\partial N(\rho_n)$ is $\partial M\setminus\lambda^{(p)}$ and the push forward of the path metric defines a hyperbolic metric $s_n$ on $\partial M\setminus\lambda^{(p)}$. We are going to show that these metrics converge.

 \begin{lemma}	\label{conval}
 	There is a multi-curve $L^0\supset \lambda^{(p)}$ such that, up to extracting a subsequence, $\ell_{s_n}(L^0)\longrightarrow 0$ and for any component $S$ of $\partial M\setminus  L^0$ the restriction of $(s_n)$ to $S$ converges to a complete hyperbolic metric.
 \end{lemma}  
 
\begin{proof}
Consider a pants decomposition $P$ given by Lemma \ref{pantssmallmeasure}, such that $i(P,\breve c(\dot\lambda_\infty))\leq\eps$ for some small $\eps>0$ and a transverse $c_\mu$ for every exceptional minimal $\mu$ of $\breve c(\dot\lambda_\infty)$ also with $i(c_\mu,\breve c(\dot\lambda_\infty))\leq\eps$. From the Small Bending Lemma \cite[Lemme A.1]{lecuire:plissage} we deduce the following : 

\begin{claim}	\label{claim:bounded}
	The sequence $(\ell_{s_n}(c))$ is bounded for any leaf $c$ of $P$ and any transverse $c=c_\mu$.
\end{claim}

\begin{proof}
By \cite[Lemme A.1]{lecuire:plissage}, there is $K$ such that $\ell_{s_n}(c)\leq K(\ell_{\rho_n}(c)+1)$. The conclusion follows from the convergence of $(\rho_n)$.
\end{proof}

As a first consequence, we deduce from the collar Lemma that if $c_n$ is a sequence of simple closed curves on $\partial M$ such that $\ell_{s_n}(c_n)\longrightarrow 0$ then $c_n$ is eventually a leaf of $P$ and is disjoint from any transverse. We can now define $L^0\subset P$ as the maximal multi-curve such that $\ell_{s_n}(c)\longrightarrow 0$ for any leaf $c$ of $L^0$, with the convention that $\ell_{s_n}(c)=0$ for $c\subset\lambda^{(p)}$.

To simplify the notation, from now on, the restriction of $P$ to $S$ will also be denoted by $P$. Extract a subsequence such that $(\ell_{s_n}(c))$ converges for any leaf $c$ of $P$. Then the restriction of $s_n$ to each component of $S\setminus  P$ converges. It follows that there is a sequence of homeomorphisms $\phi_n:S\to S$ which are products of Dehn twists along the leaves of $P$ and a subsequence such that $(\phi_n^*s_n^S)$ converges. If we denote by $D_c:S\to S$ the right Dehn twist around $c$, then for each leaf $c$ of $P$, there is $p_n(c)\in\Z$ such that $\phi_n$ is the product of $D_c^{p_n(c)}$ over the leaves of $P$.

Let $c$ be a leaf of $P$. If $i(c,\lambda)>0$, then $c$ intersects a minimal sublamination $\mu$ of $\lambda$. Lemma \ref{pantssmallmeasure} and Claim \ref{claim:bounded} provides us with a transverse $c_\mu$ such that $\ell_{s_n}(c_\mu)$ is bounded. It follows that $p_n(c)$ is bounded.

We may now assume that $i(c,\lambda)=0$. It follows that $i(\lambda_n,c)\longrightarrow0$ and by \cite[Lemme A.1]{lecuire:plissage}, $\ell_{s_n}(c)\longrightarrow\ell_{\rho_\infty} (c)$. In particular $(\ell_{s_n}(c))$ converges and since $S$ is a component of $\partial M\setminus  L^0$, $\lim(\ell_{s_n}(c))>0$ and $\rho_\infty(c)$ is not parabolic.

Let $\kappa\subset \partial M$ be a simple arc such that $\partial\kappa\subset P$, that $\kappa$ intersects $c$ in a single point $x$, that $\mathring{\kappa}\cap P=\{x\}$ and that $\kappa\setminus  P$ is disjoint from $\lambda$. Let $c_1$ and $c_2$ be the leaves of $P$ that contains the endpoints of $\kappa$ (we may have $c_1=c_2$) and denote by $\kappa_1$ and $\kappa_2$ the closures of the components of $\mathring{\kappa}\setminus  x$. Consider the loops $a=\kappa_1c_1\kappa_1^{-1}$ and $b=\kappa_2c_2\kappa_2^{-1}$. We are in the situation of figure \ref{matsu}. 

\begin{figure}[hbtp]
	\psfrag{a}{$a$}
	\psfrag{b}{$b$}
	\psfrag{x}{$x$}
	\psfrag{c}{$c$}
	\psfrag{}{$F^1$}
	\psfrag{}{$F^2$}
	\centerline{\includegraphics{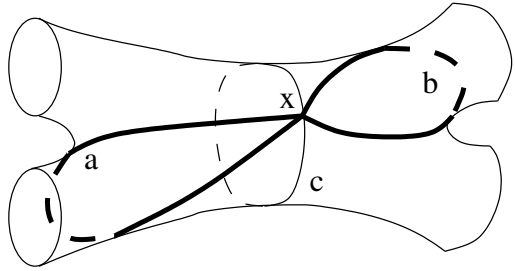}}
	\caption{Controlling the number of Dehn twists}
	\label{matsu}
\end{figure}

Let us choose $\eps$ so that $\eps<\frac{\eta}{2}\leq 2\pi$ where $\eta$ comes from the assumption that $\lambda\in\mathcal{P}(M)$ (condition $b)$ in definition \ref{def:pm}). Since $i(\lambda,P)\leq\eps$, it follows from condition $c)$ of definition \ref{def:pm} that $a$, $b$ and $c$ are not trivial and from condition $b)$ that they are primitive (i.e. if there is $d\in \pi_1(M)$ such that $d^p=a$, $b$ or $c$, then $p=1$): if a simple closed curve $e\subset\partial M$ is not primitive, then, by the Annulus Theorem  (see \cite[Theorem IV.3.1]{jaco:shalen:seifert}), two disjoint copies of $e$ bounds an essential annulus. Seeking a contradiction, let us assume that $p_n(c)\longrightarrow\infty$.

Since $(\phi_{n}^* s_n^S)$ converges, $(\rho_n(\phi_n(ab))=\rho_n(ac^{p_n(c)}b c^{-p_n(c)}))$ converges and, since $(\rho_n(a))$ converges, so does $(\rho_n(c^{p_n(c)}b c^{-p_n(c)}))$. Let $\ell_n\subset\Hp^3$ be the axis of $\rho_n(c)$ and let $x_n\in\ell_n$ be a sequence converging to $x_\infty\in\ell_\infty$. Since $\rho_n(c)$ tends to $\rho_\infty(c)$ which is hyperbolic, if $p_n(c)\longrightarrow\infty$, then $\rho_n(c)^{-p_n(c)}(x_n)$ tends to a fixed point of $\rho_\infty(c)$. On the other hand,  $d_{\Hp^3}(\rho_n(c)^{-p_n(c)}(x_n),\rho_n(b)\rho_n(c)^{-p_n(c)}(x_n))=d_{\Hp^3}(x_n,\rho_n(c)^{p_n(c)}\rho_n(b)\rho_n(c)^{-p_n(c)}(x_n))$ is bounded. This implies that $\rho_n(c)^{-p_n(c)}(x_n)$ tends to a fixed point of $\rho_\infty(b)$ and that $\rho_\infty(b)$ and $\rho_\infty(c)$ have a common fixed point. Since $\rho_\infty(\pi_1(M))$ is a discrete group, $\rho_\infty(b)$ and $\rho_\infty(c)$ commute. Since $\rho_\infty(c)$ is hyperbolic and $b$ and $c$ are primitive, $\rho_\infty(b)=\rho_\infty(c)$. It follows from \cite{waldhausen:deutsch} that $c$ and $c_2$ bound an essential annulus $E$. Since $i(\breve c(\dot\lambda_\infty),\partial E)\leq i(\breve c(\dot\lambda_\infty),P)\leq\eta$, this contradicts condition $b)$.
\end{proof}

\subsection{Convergence of the boundary of the convex core}
We continue working under the assumptions that were made at the beginning of \textsection 4, namely $(\dot\lambda_n)$ converges to $\dot\lambda_\infty$ and $(\rho_n)$ converges to $\rho_\infty$. Using the convergence of the metrics on the boundary, we are going to construct a convex pleated surface $f_\infty:S\rightarrow\Hp^3/\rho_\infty(\pi_1(M))$. 

As observed earlier, given an embedding $e_n:N(\rho_n)\to M$ associated to $\sigma_n$ the image of $\partial N(\rho_n)$ is    $\partial M\setminus\lambda^{(p)}$. The inverse of the restriction of $e_n$ to $\partial N$, $f_n=e_n^{-1}:\partial M\setminus \lambda^{(p)}\rightarrow \partial N(\rho_n)\subset\Hp^3/\rho_n(\pi_1(M))$, viewed as a map to $\Hp^3/\rho_n(\pi_1(M))$ is a convex pleated surface whose bending measured geodesic lamination is $\lambda_n\setminus \lambda^{(p)}$. We will show that this sequence of pleated surfaces $(f_n)$ tends to a pleated surface $f_\infty:\partial M\setminus L^0\to N(\rho_\infty)$, where $L^0\supset \lambda_\infty^{(p)}$ is the multi-curve defined in Lemma \ref{conval}. 

\begin{lemma}	\label{surpli}
Let $S$ be a connected component of $\partial M\setminus L^0$ and $S'$ the component of $\partial M\setminus \lambda^{(p)}$ containing $S$. Up to extracting a subsequence, there is a sequence of base points $z_n\subset S$ such that the pleated surface $f_n:S\rightarrow\Hp^3/\rho_n(\pi_1(M))$ converge to a convex pleated surface $f_\infty:S\rightarrow\Hp^3/\rho_\infty(\pi_1(M))$ which is homotopic to the inclusion map. Furthermore $f_\infty(S)$ is a component of $\partial N(\rho_\infty)$.
\end{lemma}

\begin{proof}
To extract from $(f_n)$ a converging subsequence, we will follow the ideas of \cite{canary:epstein:green}. For this, we need to show that the $f_n(S)$ all intersect the image in $\Hp^3/\rho_n(\pi_1(M))$ of a compact subset of $\Hp^3$. Let $x_n$ be the image in $\Hp^3/\rho_n(\pi_1(M))$ of the base point $O\in\Hp^3$. With a classical argument using simplicial hyperbolic surfaces, Bonahon-Otal proved :

\begin{lemma}[\cite{bonahon:otal}, lemma 17]	\label{kaenne}
There is a sequence of arcs $k_n$ joining $f_n(S)$ to the base point $x_n$ with uniformly bounded length.	\hfill $\Box$
\end{lemma}

\indent
Let $y_n=\partial k_n\setminus  x_n$ and $z_n=f_n^{-1}(y_n)$. Using Lemmas \ref{conval} and \ref{kaenne} and the arguments of \cite[\S 5.2]{canary:epstein:green}, we get that, up to extracting a subsequence, the sequence of pleated surfaces $f_n:S'\rightarrow \Hp^3/\rho_n(\pi_1(M))$ with basepoint $z_n$ converges to a pleated surface $f_\infty:S\rightarrow \Hp^3/\rho_\infty(\pi_1(M))$. Since all the $f_n$ are homotopic to the inclusion map, their limit $f_\infty$ is also homotopic to the inclusion map.

By \cite[Lemme 21]{bonahon:otal} (see also \cite[Lemma 3.1]{lecuire:continuity}), either $f_\infty$ is a covering onto its image and there is a convex set $C_{f_\infty}\subset\Hp^3/\rho_\infty(\pi_1(M))$ whose boundary is $f_\infty(S)$, or $f_\infty(S)$ lies in a totally geodesic surface  $\Hp^3/\rho_\infty(\pi_1(M))$. In the latter case $\rho_\infty(\pi_1(M))$ is a Fuchsian group and the restriction of $f_\infty$ to $S\setminus |\dot\lambda_\infty|$ is a two-sheeted covering of the interior of $f_\infty(S)$ considered as a two dimensional surface with boundary (see \cite[Lemma 3.8]{lecuire:continuity} for more details). This is only possible if $M$ is an interval bundle over a closed surface $I\times F$ and $|\dot\lambda_\infty|$ is a section of the bundle $\partial S\times I$. This would contradict the fact that $\dot\lambda_\infty$ lies in $\mathcal{P}(M)/\mathcal{R}$. Therefore $f_\infty$ is a covering onto its image and it follows from \cite[Lemma 4.4]{lecuire:continuity} that $f_\infty$ is a homeomorphism into its image (see also \cite[Lemme 3.5.2]{lecuire:thesis}).\\
\indent
Since $f_\infty$ is a pleated surface, we have $f_\infty(S)\subset N(\rho_\infty)$. On the other hand, the preimage $\tilde C_{f_\infty}\subset\Hp^3$ of $C_{f_\infty}$ under the covering projection $\Hp^3\to \Hp^3/\rho_\infty(\pi_1(M))$ is a $\rho_\infty(\pi_1(M))$-equivariant convex set. It follows from the minimality of the preimage $\tilde N(\rho_\infty)$ of $N(\rho_\infty)$ that $\tilde N(\rho_\infty)\subset \tilde C_{f_\infty}$. We have thus proved that the surface $f_\infty(S)$ is a component of the boundary of $N(\rho_\infty)$.
\end{proof}

We construct $f_\infty$ for each component $S$ of $\partial M\setminus  L_0$ and denote by $f_\infty:\partial M\setminus  L_0\rightarrow \partial N(\rho_\infty)$ the resulting map. By the above, the image of $f_\infty$ is a union of components of $\partial N(\rho_\infty)$ and $f_\infty$ is a local homeomorphism. Let us show that it is a global homeomorphism

\begin{claim}	\label{hom}
The map $f_\infty: \partial M\setminus  L^0 \rightarrow \partial N(\rho_\infty)$ is injective.
\end{claim}
\begin{proof}
Assume the contrary.  We have shown above that the restriction of $f_\infty$ to any component of $\partial_{\chi<0} M\setminus L_0$ is a homeomorphism into its image. Hence there are two connected components $S$ and $S'$ of $\partial_{\chi<0} M\setminus L_0$ such that $f_\infty(S)\cap f_\infty (S')\neq\emptyset$. Since $f_\infty(S)$ and $f_\infty(S')$ are connected components of $\partial N(\rho_\infty)$, they intersect if and only if they are equal. Since $N(\rho_n)$ lies between $f_n(S)$ and $f_n(S')$, this is possible only if $f_\infty(S)$ and $f_\infty(S')$ are totally geodesic and if $\rho_\infty(\pi_1(M))$ is Fuchsian (compare with \cite[Lemma 21]{bonahon:otal}). In this case, $M$ is homeomorphic to $S\times I$ and $|\dot\lambda_\infty|$ lies in a section of the bundle $\partial S\times I$ contradicting again the fact that $\dot\lambda_\infty$ lies in $\mathcal{P}(M)/\mathcal{R}$.
\end{proof}

\subsection{Limit metric and quasi-isometries}

We want to work with the subsequence produced by Lemma \ref{surpli} so we will now assume that $(\dot\lambda_n)$ converges to $\dot\lambda_\infty$, that $(\rho_n)$ converges to $\rho_\infty$ and that $(f_n)$ converges to $f_\infty$. Following \cite{bonahon:otal}, we prove that $\rho_\infty$ is geometrically finite and is associated to a complete hyperbolic metric on $int(M)$.

\begin{lemma}	\label{lm:limit:geom:finite}
The limit representation $\rho_\infty$ is geometrically finite and associated to a complete hyperbolic metric on $\IM$.
\end{lemma}

\begin{proof}
In this proof we will need to consider the torus components of $\partial M$. So let us temporarily go back to the usual definition of $\partial M$ which will be decomposed as $\partial M=\partial_{\chi<0} M\cup\partial_{\chi=0} M$.

To each leaf $c$ of $L^0$ correspond two cusps of $\partial_{\chi<0} M\setminus L^0$. The images of these two cusps under $f_\infty$ are two totally geodesic open annuli tending to a rank one cusp of $\Hp^3/\rho_\infty(\pi_1(M))$. By condition $b)$ for each rank one cusp, there are only two such annuli. Remove these two annuli from $f_\infty(\partial_{\chi<0} M\setminus L^0)$ and join the two boundary components of the remaining surface by a compact annulus. Doing this for each leaf of $L^0$, we get a compact surface $F_\infty\subset \Hp^3/\rho_\infty(\pi_1(M))$. In this construction, two components of $\partial_{\chi<0} M\setminus L^0$ have their images joined by an annulus if and only they are adjacent. It follows that $f_\infty$ can be turned into a homeomorphism $f:\partial_{\chi<0} M\to F_\infty$. If we denote by $F'_\infty$ the surface obtained by adding to $F_\infty$ the boundary of a neighborhood of the rank $2$ cusps of $\Hp^3/\rho_\infty(\pi_1(M))$, it is easy to extend $f$ to a homeomorphism $f:\partial M\rightarrow F'_\infty$ which is homotopic to the inclusion. Since we have a homotopy equivalence between $M$ and $\Hp^3/\rho_\infty(\pi_1(M))$ then $F'_\infty$ is homologous to $0$ and therefore bounds a compact cycle $C_\infty$. Then $f$ extends to a homotopy equivalence $f:M\to C_\infty$ and by \cite[Theorem 3.2]{evans:maps} (see also \cite[Corollary X.8]{jaco:three_manifold_topology} and \cite{waldhausen:irreducible}) we can choose this equivalence to be a homeomorphism.

 Adding some finite volume sets (which are rank $2$ cusps and ``slices'' of rank $1$ cusps) to $C_\infty$ we get a new cycle $C'_\infty$ which is bounded by $f_\infty(\partial_{\chi<0} M\setminus  L^0)$ and has finite volume. By Claim \ref{hom}, $f_\infty$ is a homeomorphism into the union of some components of $\partial N(\rho_{\infty})$. If a component $e$ of $\partial N(\rho_\infty)$ did not lie in $f_\infty(\partial_{\chi<0} M\setminus L^0)$ then the component $E$ of $\Hp^3/\rho_\infty(\pi_1(M))\setminus N(\rho_\infty)$ bounded by $e$ would lie in $C_\infty$. But, since the nearest point projection $\Hp^3/\rho_\infty(\pi_1(M))\setminus N(\rho_\infty)\longrightarrow \partial N(\rho_\infty)$ does not increase distances, $E$ has infinite volume. So $E$ does not lie in $C_\infty$. It follows that we have $f_\infty(\partial_{\chi<0} M\setminus L^0)=\partial N(\rho_\infty)$ and that $\rho_\infty$ is geometrically finite.

Since geometrically finite ends and cusps are products (\cite{marden:kg} and \cite{dbaepstein:marden}), $f$ can be extended to a homeomorphism $g:\IM\to\Hp^3/\rho_\infty(\pi_1(M))$ that defines a hyperbolic metric $\theta_\infty$ on $\IM$. By construction, $\rho_\infty$ is associated to $\theta_\infty$.
\end{proof}

With the distant goal to upgrade the convergence of representations to a convergence of (isotopy classes of) metrics we establish the existence of the desired biLipschitz maps.

\begin{lemma}	\label{lm:strong}
Let $x_n\in M_n=\Hp^3/\rho_n(\pi_1(M))$ be the projection of the origin $O\in\Hp^3$ for $n\in\N\cup\{\infty\}$. There are $k_n\longrightarrow 1$ and $r_n\longrightarrow\infty$ and homeomorphisms $u_n:M_\infty\to M_n$ such that $u_n(x_\infty)=x_n$ and that the restriction of $u_n$ to $B(x_\infty,r_n)$ is a $k_n$-bilipschitz homeomorphism onto its image.
\end{lemma}

\begin{proof}
It has been proved in \cite[Lemma 4.2]{lecuire:continuity} that there is a uniform lower bound on the length of meridians on $\partial N(\sigma_n)$, i.e. curves on $\partial N(\sigma_n)$ that are homotopically trivial in $N(\sigma_n)$ but not in $\partial N(\sigma_n)$. It follows then from \cite{canary:conformal:core} that there is a uniform lower bound on the injectivity radii of the domains of discontinuity $\Omega_{\rho_n}$ of the Kleinian groups $\rho_n(\pi_1(M))$. We deduce then from \cite{canary:poincare:core} that given a simple closed curve $c\subset\partial M$, the length of $c$ in $\Omega_{\rho_n}$ goes to $0$ when $n$ goes to $\infty$ if and only if $\ell_{s_n}(c)\longrightarrow 0$. By Lemma \ref{surpli}, $f_n(\partial M)$ converges to $\partial N(\rho_\infty)$ and $\ell_{s_n}(c)\longrightarrow 0$ if and only if $\rho_\infty(c)$ is a parabolic isometry. We then conclude from \cite{kleineidam:strong} that $\rho_n(\pi_1(M))$ converges strongly to $\rho_\infty(\pi_1(M))$, i.e. (cf. \cite[Theorem I.3.2.9 and  Corollary I.3.2.11]{canary:epstein:green}) there are sequences of numbers $k_n\longrightarrow 1$ and $r_n\longrightarrow\infty$ and maps $u_n:B(x_\infty,r_n)\to \Hp^3/\rho_n(\pi_1(M))$ such that $u_n(x_\infty)=x_n$ and that $u_n$ is a $k_n$-bilipschitz homeomorphism onto its image.

Let us denote by $M_n^{\geq\eps}$ the $\eps$-thick part of $M_n$ for $n\in\overline{\N}$. It follows from Lemmas \ref{conval} and \ref{surpli} that given $\eps$, there is a uniform bound on the diameter of $N(\rho_n)\cap M_n^{\geq\eps}$. Hence for $n$ large enough $N(\rho_n)\cap M_n^{\geq\eps}\subset u_n(B(x_\infty,r_n))$ and $N(\rho_\infty)\cap M_\infty^{\geq\eps}\subset B(x_\infty,r_n)$. Since geometrically finite ends and cusps are products (\cite{marden:kg} and \cite{dbaepstein:marden}), the homeomorphism $u_n$ is the restriction of a homeomorphism $u_n:M_\infty\to M_n$.
\end{proof}

\begin{rem}
To deduce that $\rho_n(\pi_1(M))$ converges strongly to $\rho_\infty(\pi_1(M))$, we could also have used \cite[Theorem 4.2]{jorgensen:marden} by proving that $\Omega_{\rho_n}$ converges to $\Omega_{\rho_\infty}$ in the sense of Carathéodory. This last convergence follows from the fact that $f_n(\partial M)$ converges to $\partial N(\rho_\infty)$ by tracking the behavior of support planes. The use of the work of Kleineidam was chosen because it leads to a more concise proof.
\end{rem}

\section{The action of Mod(M)}

Let $h_n:\IM\to M_n=\Hp^3/\rho_n(M)$ be a homeomorphism associated to $\sigma_n$ and $g:\IM\to M_\infty$ a homeomorphism associated to $\theta_\infty$.
To prove the convergence of $(\sigma_n)$ it remains to study the isotopy classes of the maps $h_n$ and $u_n\circ g$ ($u_n$ is the homeomorphism provided by Lemma \ref{lm:strong}). To do that we are going to study the action of the modular group on  $\mathcal{P}(M)/\mathcal{R}$. Although this may not seem to be the most efficient way to prove the convergence of $(\sigma_n)$, combined with Lemma \ref{lm:strong} and the continuity of the bending map (\cite{lecuire:continuity}), it leads to a short proof of the properness of the bending map, and it also allows us to understand the dynamic of the action of the $Mod(M)$ on $\mathcal{D}(M)$.

The {\em modular group} $Mod (M)$ is the group of isotopy class of homeomorphism $M\to M$. The restriction of each isotopy class to $\partial M$ defines an embedding of $Mod(M)$ into the mapping class group $Mod(\partial M)$. Since measured geodesic laminations are defined up to isotopy, $Mod(\partial M)$ (and hence $Mod(M)$) acts on $\ML(\partial S)$ in the obvious way: if $\phi:\partial M\to\partial M$ is a homeomorphism and $[\phi]$ denotes its mapping class then $[\phi]\lambda=\phi(\lambda)$ for any $\lambda\in\ML(\partial M)$. It is easy to see that this action projects to an action of $Mod(\partial M)$ on $\ML/\mathcal{R}$.

 We will now explain how the result of the previous section can be used to study the action of $Mod(M)$ on $\mathcal{D}(\partial M)$

\begin{lemma}	\label{lm:action}
Let $M$ be a compact orientable hyperbolic $3$-manifold, the action of $Mod (M)$ on $\mathcal{P}(M)/\mathcal{R}$ is properly discontinuous.
\end{lemma}

\begin{proof}
Consider $(\dot\lambda_n)\in \mathcal{P}(M)/\mathcal{R}$ and $(\phi_n)\in Mod(M)$ such that $(\dot\lambda_n)$ converges to $\dot\lambda_\infty\in\mathcal{P}(M)/\mathcal{R}$ and that $(\phi_n(\dot\lambda_n))$ converges to $\mu_\infty\in \mathcal{P}(M)/\mathcal{R}$. We are going to show that the homeomorphisms $\phi_n$ lie in finitely many isotopy classes.

We start with a subsequence that we still denote with the same index $_n$ and we will not change the index when extracting further subsequences. Since $\dot\lambda_n\in \mathcal{P}(M)/\mathcal{R}$, $\breve c(\dot\lambda_n)$ is the representative lying in $\mathcal{P}(M)$. By \cite{lecuire:plissage}, for any $n\in\N$, there is a geometrically finite metric $\sigma_n\in\mathcal{GF}(M)$ whose bending measured lamination is $\breve c(\dot\lambda_n)$.
 
By Proposition \ref{algebrique}, up to extracting a subsequence, there are two sequences representations associated to $\sigma_n$ and $\phi_*(\sigma_n)$ respectively that converge. By Lemma \ref{conval}, there is a multi-curve $L^0$, resp. $N^0$ such that, up to extracting a further subsequence, the restrictions of $(s_n)$, resp. $(\phi_{n*}(s_n))$, to each component of $\partial M\setminus  L^0$, resp. $\partial M\setminus N^0$,  converge to complete hyperbolic metrics.

In particular, $\phi_n(L^0)=N^0$ for $n$ large enough.

Consider a weighted multi-curve $\gamma$ whose support contains $L^0$. By choosing $\gamma$ close to $\breve c(\dot\lambda_\infty)$, we may assume that $\gamma\in\mathcal{P}(M)$. Since the restrictions of $(s_n)$, resp $(\phi_{n*}(s_n))$, to each component of $\partial M\setminus  L^0$, resp. $\partial M\setminus  N^0$, converge, $(\ell_{s_n}(\gamma))$ and $(\ell_{\phi_{n*}(s_n)}(\gamma))$ are bounded. It follows that $\phi_n(\gamma)$ is finite up to isotopy, and up to extracting a further subsequence, $\phi_n(\gamma)$ is isotopic to $\phi_m(\gamma)$ for any $n,m\in\N$.

Let $P$ be a regular neighborhood of the support of $\gamma$ and let $F$ be the closure of $\partial M\setminus P$. The components of $P$ and $F$ form a finite set of compact $2$-manifolds whose intersections are either empty or $1$-manifolds. Johannson called such a set a boundary-pattern  (see \cite{johannson:jsj}). Following Johannson, we will denote this boundary-pattern by $\underline{\underline{m}}$. It is complete since any point of $\partial M$ lies in some element of $\underline{\underline{m}}$. 
By properties $a)$ and $b)$, any essential disc has at least four intersections with the boundaries of the elements of $\underline{\underline{m}}$, Joannson calls a boundary-pattern with this property useful. By property $c)$, no essential annulus can be disjoint from the boundaries of the elements of $\underline{\underline{m}}$, such a boundary-pattern is said to be simple. From the previous paragraph we get that $\phi_n\circ\phi_m^{-1}(\underline{\underline{m}})$ is isotopic to $\underline{\underline{m}}$. It follows then from \cite[Proposition 27.1]{johannson:jsj} that $\{\phi_n\circ\phi_m^{-1}, n,m\in\N\}$ is finite up to isotopy, hence the homeomorphisms $\phi_n$ lie in finitely many isotopy classes. 

We have shown that any subsequence contains a further subsequence such that homeomorphisms $\phi_n$ lie in finitely many isotopy classes. This is only possible if the original sequence contains only finitely many classes of homeomorphisms.  Thus we have shown that the action of $Mod(M)$ on $\mathcal{P}(M)/\mathcal{R}$ is properly discontinuous.
\end{proof}

\setcounter{section}{1}
\setcounter{theorem}{1}

Lemma \ref{lm:action} allows us to conclude the proof of Theorem \ref{two}

\begin{theorem}
	The map $b_\mathcal{R}$ from $\mathcal{GF}(M)$ to $\mathcal{P}(M)/\mathcal{R}$ is proper.
\end{theorem}

\begin{proof}
	Let $\sigma_n$ be a sequence of geometrically finite metrics on $\IM$ with bending measured laminations $\lambda_n$. Assume that the projections $\dot\lambda_n$ on $\mathcal{P}(M)/\mathcal{R}$ converges to $\dot\lambda_\infty\in\mathcal{P}(M)/\mathcal{R}$. Up to extracting a subsequence, by Proposition \ref{algebrique} we may chose representations $\rho_n$ associated to $\sigma_n$ such that $(\rho_n(g))$ converges to some $\rho_\infty(g)$ for any $g\in\pi_1(M)$. By Lemma \ref{lm:limit:geom:finite} $\rho_\infty(M)$ is geometrically finite and there is a homeomorphism $g:\IM:\to M_\infty=\Hp^3/\rho_\infty(M)$ defining a metric $\theta_\infty$ on $M$. By Lemma \ref{lm:strong} there are homeomorphisms $u_n:M_\infty\to M_n$ and sequences $k_n\longrightarrow 1$ and $r_n\longrightarrow\infty$ such that $u_n(x_\infty)=x_n$ and that the restriction of $u_n$ to $B(x_\infty,r_n)$ is a $k_n$-bilipschitz homeomorphism onto its image. 
	
	Let $h_n:\IM\to M_n=\Hp^3/\rho_n(M)$ be a homeomorphism associated to $\sigma_n$. The homeomorphism $g_n=u_n\circ g:\IM\to M_n$ defines a metric $\theta_n$ and by the previous paragraph $\theta_n$ converges to $\theta_\infty$. Set $\phi_n=g_n^{-1}\circ h_n:\IM\to\IM$ and extend it to a homeomorphism $\phi_n:M\to M$. By construction $\theta_n=\phi_{n*}\sigma_n$ and $\phi_n(\lambda_n)$ is the bending lamination of $\theta_n$. It follows from \cite{lecuire:continuity} that $(\phi_n(\dot\lambda_n))$ converges to $\dot\mu_\infty$ where $\mu_\infty\in\mathcal{P}(M)$ is the bending lamination of $\theta_\infty$. By Lemma \ref{lm:action}, the homeomorphisms $\phi_n$ lie in finitely many isotopy classes. Since $\theta_n=\phi_{n*}\sigma_n$ converge, it follows that $(\sigma_n)$ has a converging subsequence.
\end{proof}

We will conclude this article with the proof of Theorem \ref{tree}.

\begin{theorem}
When $M$ is not a genus two handlebody, the action of $Mod (M)$ on $\mathcal{D}(M)$ is properly discontinuous.
\end{theorem}
 
\begin{proof}
We will use Lemma \ref{lm:action} and we start the proof in the same manner. Consider $(\lambda_n)\subset \mathcal{D}(M)$ and $(\phi_n)\subset Mod(M)$ such that $(\lambda_n)$ converges to $\lambda_\infty\in\mathcal{D}(M)$ and that $(\phi_n(\lambda_n))$ converges to $\mu_\infty\in \mathcal{D}(M)$.
	
Since $\lambda_\infty\in\mathcal{D}(M)$ and $\mathcal{D}(M)$ is open (\cite[Lemma 4.1]{lecuire:masur}), there is $\eta>0$ such that $i(\lambda_n,\partial D)>\eta$ for any essential disc $D$ and for $n$ large enough. Let $\frac{2\pi}{\eta} \lambda_n$ be the measured geodesic lamination obtained by rescaling the measure of  $\lambda$ by $\frac{2\pi}{\eta}$. It is shown in the proof of \cite[Lemma 3.5]{lecuire:masur} that $c(\frac{2\pi}{\eta} \lambda_n)\in\mathcal{P}(M)$. Since $(\lambda_n)$, resp. $(\phi_n(\lambda_n))$, converges to $\lambda_\infty\in\mathcal{D}(M)$, resp. $\mu_\infty\in \mathcal{D}(M)$, the projections of $c(\frac{2\pi}{\eta} \lambda_n)$, resp. $\phi_n(c(\frac{2\pi}{\eta}\lambda_n))$, to $\mathcal{P}(M)/\mathcal{R}$
converges to the projection of $c(\frac{2\pi}{\eta} \lambda_\infty)$, resp. $c(\frac{2\pi}{\eta} \mu_\infty)$.

By Lemma \ref{lm:action}, the action of  $Mod (M)$ on $\mathcal{P}(M)$ is properly discontinuous. Hence the $\phi_n$ lie in finitely many isotopy classes.
\end{proof}

\bibliographystyle{math}
\bibliography{refplus}

\end{document}